\documentclass[10pt]{amsart}
\usepackage{amsthm,amssymb,amsmath,amsfonts,mathrsfs,amscd,yfonts,stmaryrd}
\usepackage{url,tikz,enumitem,caption,stmaryrd,comment,mathrsfs,marginnote,colonequals}
\usepackage[utf8]{inputenc}
\usepackage[T1]{fontenc}
\usepackage{libertine}
\usetikzlibrary{matrix,positioning,arrows,shapes,chains,calc,automata}
\usetikzlibrary{decorations.pathmorphing}
\usetikzlibrary{arrows,automata,positioning,knots}
\definecolor{darkgreen}{rgb}{0,0.5,0}
\usepackage[pagebackref=true,
        colorlinks, citecolor=darkgreen,
]{hyperref}

\newlength{\proofmargin}
\renewcommand*{\backref}[1]{}
\renewcommand*{\backrefalt}[4]{%
  \ifcase #1 %
    \relax
  \or
    $\uparrow$#2.%
  \else
    $\uparrow$#2.%
  \fi%
}

  \renewenvironment{thebibliography}[1]{%
    \begin{oldthebibliography}{#1}%
      \setlength{\parskip}{0ex}%
      \setlength{\itemsep}{0.1ex}%
      \setlength{\labelwidth}{1.4cm} 
  }%
  {%
    \end{oldthebibliography}%
  }

\DeclareMathOperator{\Q}{\mathbf{Q}}
\DeclareMathOperator{\F}{\mathbf{F}}
\DeclareMathOperator{\Z}{\mathbf{Z}}
\DeclareMathOperator{\G}{\mathbf{G}}

\DeclareMathOperator{\C}{\mathbf{C}}
\DeclareMathOperator{\PP}{\mathbf{P}}

\DeclareMathOperator{\Fil}{\mathrm{Fil}}

\DeclareMathOperator{\cX}{\mathcal{X}}

\DeclareMathOperator{\cH}{\mathcal{H}}
\DeclareMathOperator{\cO}{\mathcal{O}}

\DeclareMathOperator{\GL}{\mathrm{GL}}

\DeclareMathOperator{\D}{\mathrm{D}}
\DeclareMathOperator{\U}{\mathrm{U}}
\DeclareMathOperator{\Fr}{\mathrm{Fr}}

\DeclareMathOperator{\HH}{\mathrm{H}}
\DeclareMathOperator{\Hf}{\mathrm{H}_f^1}
\DeclareMathOperator{\AJb}{\mathrm{AJ}_b}

\DeclareMathOperator{\cA}{\mathcal{A}}

\DeclareMathOperator{\cC}{\mathcal{C}}

\DeclareMathOperator{\lra}{\longrightarrow}

\DeclareMathOperator{\res}{\mbox{\scriptsize res}}

\DeclareMathOperator{\Tr}{\mathrm{Tr}}

\DeclareMathOperator{\Ab}{\mbox{\scriptsize ab}}
\DeclareMathOperator{\Gal}{\mbox{\scriptsize Gal}}

\DeclareMathOperator{\an}{\mbox{\scriptsize an}}
\DeclareMathOperator{\et}{\mbox{\scriptsize \'et}}
\DeclareMathOperator{\cris}{\mbox{\scriptsize cris}}
\DeclareMathOperator{\dR}{\mbox{\scriptsize dR}}
\DeclareMathOperator{\GM}{\mbox{\scriptsize GM}}

\newcommand{\per}{\mathrm{per}_p}

\newtheorem{theorem}{Theorem}[section]
\newtheorem{lemma}[theorem]{Lemma}

\newtheorem{definition}[theorem]{Definition}
\newtheorem{conjecture}[theorem]{Conjecture}

\theoremstyle{remark}
\newtheorem{algorithm}[theorem]{Algorithm}
\newtheorem{remark}[theorem]{Remark}

\newcommand{\Affq}{\operatorname{Aff}(q)}

\newcommand{\MCG}{\mathrm{MCG}}

\makeatletter
\renewenvironment{proof}[1][\proofname]%
{%
\par\pushQED{\qed}\normalfont\topsep6\p@\@plus6\p@\relax%
\begin{list}{}{\rightmargin=8pt\leftmargin=\proofmargin}%
  \item[\hskip\labelsep\bfseries#1\@addpunct{.}]\ignorespaces
}{%
\popQED\end{list}\@endpefalse%
}%
\makeatother

\title[Two recent $p$-adic approaches towards the (effective) Mordell conjecture]{Two recent $p$-adic approaches towards\\ the (effective) Mordell conjecture}

\author[Balakrishnan]{Jennifer S. Balakrishnan}
\address{\hspace{-.2in}J. S. Balakrishnan, Department of Mathematics \& Statistics, Boston University, 111 Cummington Mall, Boston, MA 02215,~USA}
\email{jbala@bu.edu}

\author[Best]{Alex J. Best}
\address{\hspace{-.2in}A. J. Best, Department of Mathematics \& Statistics, Boston University, 111 Cummington Mall, Boston, MA 02215, USA}
\email{alex.j.best@gmail.com}

\author[Bianchi]{Francesca Bianchi}
\address{\hspace{-.2in}F. Bianchi,  Bernoulli Institute,   University of Groningen,  Nijenborgh 9,  9747 AG Groningen,  The Netherlands}
\email{francescabianchi25@gmail.com}

\author[Lawrence]{Brian Lawrence}
\address{\hspace{-.2in}B. Lawrence, Department of Mathematics, University of Chicago, 5734 S University Ave, Chicago, IL 60637, USA}
\email{brianrl@math.uchicago.edu}

\author[M\"uller]{J. Steffen M\"uller}
\address{\hspace{-.2in}J. S. M\"uller,  Bernoulli Institute,   University of Groningen,  Nijenborgh 9,  9747 AG Groningen,  The Netherlands}
\email{steffen.muller@rug.nl}

\author[Triantafillou]{Nicholas Triantafillou}
\address{\hspace{-.2in}N. Triantafillou, Department of Mathematics, University of Georgia, Athens, GA 30602, USA}
\email{nicholas.triantafillou@uga.edu}

\author[Vonk]{Jan Vonk}
\address{\hspace{-.2in}J. Vonk, Institute for Advanced Study, 1 Einstein Drive, Princeton, NJ 08540, USA}
\email{vonk@ias.edu}

\begin{document}

\maketitle
\tableofcontents

\section{Introduction}
\label{sec:Intro}
\subsection{The Mordell conjecture. }Many important developments in arithmetic geometry
were motivated by the Mordell conjecture, stated nearly a century ago. Let $X$ be a smooth
projective curve, defined over the field of rational numbers $\Q$. Its set of rational
points $X(\Q)$, which consists of all the projective solutions with rational coordinates
to a finite set of equations defining $X$ in some projective space, is an interesting arithmetic quantity. In 1922, Mordell \cite{Mor22} made the following conjecture:
\begin{conjecture}[Mordell]
Suppose that $X$ is of genus at least two. Then $X(\Q)$ is finite. 
\end{conjecture}
In a monumental paper, Faltings \cite{Fal83} proved this conjecture. The method of Faltings is ingenious, and merits a thorough treatment on its own. Indeed, many such are available in the literature, see for instance \cite{CS86} for an early account. In this paper, we wish to give an introductory account of two recent alternative approaches towards this conjecture, due to Lawrence--Venkatesh \cite{LV18} and Kim \cite{Kim05, Kim09, Kim10b}. The latter method, which is usually called the method of \textit{Chabauty--Kim} or \textit{non-abelian Chabauty} in the literature, has the advantage that in some cases it has been turned into an effective method to \textit{determine} the set of rational points $X(\Q)$, and we illustrate this by presenting three new examples of modular curves where this set can be determined, due to Best, Bianchi, and Triantafillou. 

\begin{remark}Mordell's conjecture, as well as many of the results discussed below, admit
  analogues where $X$ is replaced by a smooth \textit{hyperbolic curve}, including also
  the cases of a punctured elliptic curve and $\PP^1 \backslash \{ 0,1,\infty \}$, when
  the set of rational points $X(\Q)$ is replaced by the set of $S$-integral points, where
  $S$ is a finite set of primes. In this setting, the finiteness of $S$-integral points is
  known as Siegel's theorem.
Both methods presented here are expected to apply to non-proper hyperbolic curves as well.
We discuss the $S$-unit equation in the context of \cite{LV18} below. 
  Kim~\cite{Kim05} proved the finiteness of integral points on $\PP^1 \backslash \{
    0,1,\infty \}$,  and explicit Chabauty-Kim methods for $S$-unit equations are due to Dan-Cohen and Wewers~\cite{DCW15,DCW16,
  DC}. Chabauty-Kim theory for integral points on punctured
  elliptic curves of rank 0 and 1 is discussed in~\cite{Kim10b} and~\cite{BDCKW18}.
\end{remark}

\begin{remark}For the purpose of exposition, we only consider the base field $\Q$. It should be noted that many results admit appropriate generalizations to number fields \cite{siksek:nf, D19, BBBM}. The only exception is our discussion of the method of Lawrence--Venkatesh, where field extensions play an essential role. \end{remark}

\subsection{Two recent approaches. }
After Faltings' proof, two notable new methods for proving finiteness of $X(\Q)$ for $X$ of genus $g \ge 2$ have emerged. In broad strokes, they follow a similar strategy: We start by choosing a prime $p$ at which the curve $X$ has good reduction, and we study the set of rational points through the inclusion
\begin{equation}
X(\Q) \  \subset \ X(\Q_p).
\end{equation}
For any field $K$, we write $G_K = \mathrm{Gal}(\overline{K}/K)$ for its absolute Galois group. The starting point of both the methods of Chabauty--Kim and Lawrence--Venkatesh is the association of a certain finite-dimensional Galois representation over $\Q_p$ to every point on the curve, giving maps
\begin{equation}
\rho: X(K) \ \lra \ \mathrm{Rep}(G_K),  
\end{equation}
for $K$ equal to $\Q$ or $\Q_p$. In both the approaches of Lawrence--Venkatesh and Chabauty--Kim, finiteness of the set $X(\Q)$ is obtained from the consideration of a commutative diagram of the following shape: 
\begin{equation}\label{eqn:diagram-ab}
\begin{tikzpicture}[->,>=stealth',baseline=(current  bounding  box.center)]
 \node[] (X) {$X(\Q)$};
 \node[right of=X, node distance=3.2cm]  (Xp) {$X(\Q_p)$};   
 \node[below of=X, node distance=1.5cm]  (RepG) {$\mathrm{Rep}(G_{\Q})$};
 \node[right of=RepG,node distance=3.2cm] (RepGp) {$\mathrm{Rep}(G_{\Q_p})$};
 \node[right of=RepGp,node distance=2.7cm](Fil) {$\mathrm{MF}^{\phi}/ \simeq$.};

 \path (X)  edge node[left]{\footnotesize $\rho$} (RepG);
 \path (Xp) edge node[left]{\footnotesize $\rho$} (RepGp);
 \path (X)  edge (Xp);
 \path (RepG) edge node[above]{\footnotesize $\mathrm{res}_p$} (RepGp);
 \path (RepGp) edge node[above]{\footnotesize $\D_{\cris}$}(Fil);
 \path (Xp)  edge[dashed, bend left=25] node[above right]{\footnotesize $\per$} (Fil);
\end{tikzpicture}
\end{equation}
While the nature of $\rho$ is very different in the two approaches, the horizontal maps are the same. First of all, the map from $X(\Q)$ to $X(\Q_p)$ is simply the natural inclusion, and $\mathrm{res}_p$ is the restriction of Galois representations, making the diagram commutative in both approaches. The map $\D_{\cris}$ is defined using $p$-adic Hodge theory. More precisely, it is Fontaine's crystalline Dieudonn\'e functor from $p$-adic Galois representations to filtered $\phi$-modules. Finally, $\per$ is defined to be the composite of this map with $\rho$, and will be referred to as the ($p$-adic) \textit{period map}. 

\par As mentioned above, the maps $\rho$ which feature in the methods of Lawrence--Venkatesh and Chabauty--Kim are of a very different nature, and are responsible for the drastic differences between the two approaches. They may roughly be described as follows:
\begin{itemize}
\item The method of Lawrence--Venkatesh starts by considering a family of curves $\cC \longrightarrow X$. This is a so-called \textit{Parshin family}, where the fibre $\cC_x$ of a point $x$ in $X(K)$ is itself a disjoint union of finite coverings of $X$, unramified away from the point $x$. The association $\rho$ is then simply
\begin{equation*}
\rho \ \colon \ x \ \longmapsto \ \HH^1_{\et}(\overline{\cC}_{x}, \Q_p).
\end{equation*}
A lemma of Faltings can be used to show that the number of global representations in $\rho(X(\Q_p))$ is finite. The main part of the argument of Lawrence--Venkatesh is to establish that the map $\per$ is finite-to-one. The argument starts by realizing $\per$ as the quotient of the Hodge filtration map $\Phi_p : X(\Q_p) \lra \mathrm{Gr}(g,2g)$ by the Frobenius centralizer, and showing that on every residue disk
\begin{enumerate}
\item every orbit of the Frobenius centralizer has positive codimension in $\mathrm{Gr}(g,2g)$, and 
\item the image of $\Phi_p$ is Zariski dense.
\end{enumerate}
The former is established via carefully extending the base field and exploiting the semi-linearity of the Frobenius operator, whereas the latter is established using a monodromy calculation for the family $\cC$. The finiteness of $X(\Q)$ follows easily from the above commutative diagram. \\

\item In the method of Chabauty--Kim, one chooses a rational base point $b \in X(\Q)$ and obtains the association $\rho$ by considering certain well-chosen \textit{unipotent} quotients $\mathrm{U}(b)$ of the algebraic fundamental group $\pi_1^{\et}(\overline{X} ; b)$. This choice of quotient typically depends on the specifics of the curve $X$ under consideration. The association $\rho$ in the method of Chabauty--Kim is then of the form
\begin{equation*}
\rho \ \colon \ x \ \longmapsto \ \mathrm{U}(b,x)
\end{equation*}
where $\mathrm{U}(b,x)$ is obtained by twisting the unipotent quotient $\mathrm{U}(b)$ by the path torsor $\pi_1^{\et}(\overline{X} ; b,x)$. This carries the structure of a $\Q_p$-representation of $G_K$ whenever $x$ is in $X(K)$. All these Galois representations are twists of $\mathrm{U}(b)$, whose unipotence provides a certain rigidity that is crucial for arithmetic applications. More precisely, Kim shows that the image of $\rho$ is contained in a set that naturally carries the structure of an algebraic variety, which is usually referred to as a \textit{Selmer variety}, such that the map $\mathrm{res}_p$ between the global and local Selmer varieties is algebraic. 

\par This rigidity provides us with a clear strategy to prove finiteness, in the style of the classical method of Chabauty (see below). Indeed, if we can establish that 
\begin{enumerate}
\item\label{prop:1} the global Selmer variety has positive codimension in the local Selmer variety, and
\item\label{prop:2} the image of $X(\Q_p)$ is Zariski dense,
\end{enumerate}
then the intersection of the two sets (which contains the set of rational points) must be finite. Property (\ref{prop:2}) is true in great generality, whereas (\ref{prop:1}) typically requires additional information. Note the amusing similarity with the two steps in the proof of Lawrence--Venkatesh discussed above.
\end{itemize}

\par In spite of the apparent similarity of the two strategies, the different nature of the maps $\rho$ already lays bare a crucial difference: In contrast with the unconditional proof of Lawrence--Venkatesh, an additional piece of information is needed to deduce finiteness from the method of Chabauty--Kim. Typically this either takes the form of a \textit{geometric} assumption, such as having a large N\'eron--Severi rank \cite{BD18,BD19}, or the assumption of a geometric conjecture, such as the Bloch--Kato conjecture, see \cite{Kim09}. 

\subsection{Finding rational points explicitly. }\label{sec:rat_pts_expl}
At first glance, it may seem from the above comments that the conditional nature of the proof of finiteness obtained from the method of Chabauty--Kim puts the method at a significant disadvantage, especially when compared to the unconditional proof of Lawrence--Venkatesh. However, recent developments \cite{BD18,BD19,BDMTV19} have shown that in certain examples where additional geometric information is known, the method for proving finiteness can in fact be turned into a method to \textit{explicitly determine} the finite set $X(\Q)$. 

\par To explain the ideas, we briefly remind the reader of the method of Chabauty--Coleman \cite{Cha41,Col85}, of which an excellent exposition may be found in McCallum--Poonen \cite{McP12}. In this method, one chooses a rational base point $b$ in $X(\Q)$ and attaches to every other point a torsor of the $p$-adic Tate module $V$ of the Jacobian $J$.  More precisely, if $K$ is either $\Q$ or $\Q_p$, this torsor is obtained by the composition
\begin{equation}
\label{eqn:AJ-torsor}
\rho \ : \ X(K) \ \lra \ J(K) \ \lra \ \HH^1_f(G_{K},V)
\end{equation}
where the first map is the Abel--Jacobi embedding attached to the choice of base point $b$, and the second map attaches to a point $x$ in $J(K)$ the torsor of $V$ obtained from the inverse limit of the preimages of $x$ under the multiplication-by-$p^n$ map on the Jacobian, i.e. 
\begin{equation}
\label{eqn:torsor-J}
\Q_p \otimes_{\Z_p} \left( \varprojlim_{n} \, [p^n]^{-1}(x)\right).
\end{equation}
Such torsors are classified by the cohomology group $\HH^1(G_K,V)$ and satisfy certain \textit{Selmer conditions}\footnote{We are deliberately vague about these finiteness conditions here, but mention that the discussion below can be made unconditional on the finiteness of the Tate--Shafarevich group of the Jacobian. } which are denoted by the subscript $f$. This association $\rho$ is familiar from the context of the classical method of descent, used to compute the Mordell--Weil group of the Jacobian. 

\par We now obtain the commutative diagram:
\begin{equation}
\begin{tikzpicture}[->,>=stealth',baseline=(current  bounding  box.center)]
 \node[] (X) {$X(\Q)$};
 \node[right of=X, node distance=3.2cm]  (Xp) {$X(\Q_p)$};   
 \node[below of=X, node distance=1.5cm]  (RepG) {$\HH^1_f(G_{\Q},V)$};
 \node[right of=RepG,node distance=3.2cm] (RepGp) {$\HH^1_f(G_{\Q_p},V)$};
 \node[right of=RepGp,node distance=2.7cm](Fil) {$\HH^0(X, \Omega^1_X)^{\vee}$};

 \path (X)  edge node[left]{\footnotesize $\rho$} (RepG);
 \path (Xp) edge node[left]{\footnotesize $\rho$} (RepGp);
 \path (X)  edge (Xp);
 \path (RepG) edge node[above]{\footnotesize $\mathrm{res}_p$} (RepGp);
 \path (RepGp) edge node[above]{\footnotesize $\sim$}(Fil);
 \path (Xp)  edge[dashed, bend left=25] node[above right]{\footnotesize $\per$} (Fil);
\end{tikzpicture}
\end{equation} 
representing perhaps the simplest instance of the Chabauty--Kim strategy towards the Mordell conjecture discussed above, where $\U$ is taken to be the \textit{abelianization} $V$ of the fundamental group. In this situation, the relevant filtered $\phi$-modules are classified by the dual to the space of holomorphic differentials on $X$, which is of dimension $g$, and the isomorphism is provided by the Bloch--Kato logarithm. With suitable finiteness conditions $f$, the dimension of $\HH^1_f(G_{\Q},V)$ can be bounded above by the rank $r$ of the $\Q$-rational points of the Jacobian of $X$. The discussion of how to prove finiteness of $X(\Q)$ using the method of Chabauty--Kim then specializes to the classical argument of Chabauty, who deduces finiteness under the assumption that $r < g$. 

\par Going one step further, we note that the $p$-adic period map $\per$ has the following concrete description:
\begin{equation}
\per(x) = \left(\omega \longmapsto \int_b^x \omega \right)
\end{equation}
where the integration is taken in the sense of Coleman \cite{Col85}. Our ability to compute Coleman integrals \cite{BBK10, B15, BT19} often results in an explicit determination of the set $X(\Q)$. Since there already exist several excellent expositions of this method \cite{McP12}, we will simply explain the method by showing it in action for a single example.

\par\textbf{Example. }Let $X$ be the genus 3 hyperelliptic curve with minimal model\footnote{Here $X$ is the curve of absolute discriminant and conductor both equal to $60329= 23\cdot 43\cdot 61$ from the database \cite{g3hyp}.}
\begin{equation*}
 w^2 + (z^4 + z^2 + z + 1)w = -z^5 - z^2.
\end{equation*}
A search for points with small coordinates gives that
\begin{equation}
\label{eq:known_points_chabauty}
\left\{\infty^{\pm}, (-1,-2),(-1,0),(0, -1),(0, 0)\right\}\subseteq X(\Q),
\end{equation}
where $\infty^{+}=(1:0:0)$ and $\infty^{-}=(1:-1:0)$ are the points at infinity. In order to determine the full set of rational points $X(\Q)$, we apply the Chabauty--Coleman method with $p=3$; for convenience, we work with the following model for $X$:
\begin{equation*}
X^{\prime}\colon  w^2 = z^8 + 2z^6 - 2z^5 + 3z^4 + 2z^3 - z^2 + 2z + 1.
\end{equation*}
We embed $X^{\prime}$ into its Jacobian $J$ via the Abel--Jacobi map corresponding to the base point $b=(0,1)$ in $X^{\prime}(\Q)$. A computation in {\tt Magma} \cite{BCP97} shows that the Mordell--Weil rank of $J$ is equal to $1$, and the above discussion then implies that the codimension of the image of $\mathrm{res}_3$ is at least $2$. In fact, it is precisely equal to $2$: the set $\left\{\omega_i = z^i \frac{dz}{w}:0\leq i\leq 2\right\}$ is a basis for $\HH^0(X^{\prime}, \Omega^1_{X^{\prime}})$ and we have
\begin{alignat*}{3}
\mathrm{per}_3(0,-1)(\omega_0) &\equiv 3(3 + 3^3 + 2\cdot 3^4) &&\bmod{3^6}\\
\mathrm{per}_3(0,-1)(\omega_1) &\equiv 3(1 + 3^3 + 3^4)&&\bmod{3^6}\\\mathrm{per}_3(0,-1)(\omega_2) &\equiv 3(1 + 3^2 + 2\cdot 3^3 + 2\cdot 3^4)\quad &&\bmod{3^6}.
\end{alignat*}
Thus, we may choose generators $\alpha=a_0\omega_0-a_1\omega_1$ and $\beta=b_0\omega_0-b_2\omega_2$ for the $\Q_3$-vector space
\begin{equation*}
\left\{\omega\in \HH^0(X^{\prime}, \Omega^1_{X^{\prime}}) : \mathrm{res}_3(c)(\omega)=0\ \text{for all}\ c\in \HH^1_f(G_{\Q},V)\right\}
\end{equation*} 
such that
\begin{equation}
\begin{array}{llllllll}
a_0 &\equiv & 1 + 3^3 + 3^4 & \mod{3^5} & \qquad a_1 & \equiv & 3 + 3^3 + 2\cdot 3^4 & \mod{3^5}\\
b_0 & \equiv & 1 + 3^2 + 2\cdot 3^3 + 2\cdot 3^4 & \mod{3^5} & \qquad b_2 & \equiv & 3 + 3^3 + 2\cdot 3^4 & \mod{3^5}.
\end{array}
\end{equation}

By construction, we have
\begin{equation}
\label{eq:vanishing_per}
X^{\prime}(\Q)\subseteq \{x\in X^{\prime}(\Q_3) : \mathrm{per}_3(x)(\alpha) = 0 \ \text{and}\ \mathrm{per}_3(x)(\beta) = 0\}\equalscolon \mathcal{T};
\end{equation}
a computation shows that $\mathcal{T}$ contains precisely 6 points and hence that the inclusion in \eqref{eq:known_points_chabauty} is in fact an equality. Explicitly, suppose for instance that we want to compute all $x\in \mathcal{T}$ which reduce to the point $(1:1:0)$ in $X^{\prime}(\F_3)$. For $\gamma\in \{\alpha,\beta\}$ we have
\begin{equation*}
\mathrm{per}_3(x)(\gamma) = \mathrm{per}_3(1:1:0)(\gamma)+ \int_{(1:1:0)}^x\gamma = \int_{(1:1:0)}^x\gamma;
\end{equation*}
expanding in terms of the local parameter $t=z(x)^{-1}$ and formally integrating yields
\begin{alignat*}{2}
\mathrm{per}_3(x)(\alpha)&=(2\cdot 3 + 3^2)\cdot t^2 + (2\cdot 3^{-1} + 2 + 2\cdot 3+ 3^2)\cdot t^3\quad &&\bmod{(3^3,t^4)}\\
\mathrm{per}_3(x)(\beta) &= 3\cdot t + (2\cdot 3^{-1} + 1  + 3 + 2\cdot 3^2)\cdot t^3 &&\bmod{(3^3,t^4)}.
\end{alignat*}
The $i$-th coefficient of the local expansion of $\mathrm{per}_3(x)(\alpha)$ or $\mathrm{per}_3(x)(\beta)$ has valuation bounded from below by $-\mathrm{ord}_3(i)$; from Newton polygon considerations, we deduce that
\begin{itemize}
\item $\mathrm{per}_3(x)(\alpha)$ has a double zero at $t=0$, a simple zero at some $t\in \Z_3$ which satisfies $t\equiv 2\cdot 3^2 \bmod{3^3}$, and no other zero in $3\Z_3$;
\item $\mathrm{per}_3(x)(\beta)$ has a simple zero at $t=0$, two simple zeros congruent modulo $3^2$ to  $2\cdot 3$ and $3$, respectively, and no other zero in $3\Z_3$.
\end{itemize}
Therefore, the intersection of the zero sets of $\mathrm{per}_3(x)(\alpha)$ and $\mathrm{per}_3(x)(\beta)$ in the residue disk of the point $(1:1:0)$ in $X^{\prime}(\F_3)$ is precisely $\{(1:1:0)\}\subset X^{\prime}(\Q)$. 

We emphasize that neither $\alpha$ nor $\beta$ on their own would have sufficed to determine $X(\Q)$, as each of $\mathrm{per}_3(x)(\alpha)$ and $\mathrm{per}_3(x)(\beta)$ vanishes at some points $x\in X^{\prime}(\Q_3)\setminus X^{\prime}(\Q)$ which we can only compute modulo $3^n$ for a choice of $n$. More generally, for curves $X$ of genus $g$ with rank $g-1$ Jacobians, the Chabauty--Coleman method typically provides us with only one locally analytic function whose zero set $\mathcal{T}$ contains $X(\Q)$. It is then often the case that $\mathcal{T}$ contains some points that we cannot recognize as points in $X(\Q)$. In such situations, the Mordell--Weil sieve (see \S\ref{subsec:MordellWeil} for a discussion) can often be used to prove that the $p$-adic approximations of these points that we have computed cannot be approximations of points in $X(\Q)$.

\subsection{Integral points on higher-dimensional varieties. }
Both methods (Lawrence--Venkatesh and Chabauty--Kim) could be applied to the problem of finding integral points on higher-dimensional varieties as well.
To fix ideas, let $X$ be a smooth variety over $\Q$.
If $X$ has large nonabelian fundamental group, 
one can hope to construct a nontrivial $p$-adic local system on $X$.
This will attach a Galois representation to every point $x \in X(\Q)$,
giving rise to a period map $\rho$ as above.
In the higher-dimensional setting, one can no longer conclude finiteness of integral points;
rather, these methods give the weaker result that $X(\Z)$ is not dense for the $p$-adic analytic topology. 
See \cite[Section 9]{Had11} for a result of this form.
(In dimension one, non-density for the analytic topology is equivalent to finiteness.)

It is sometimes possible to strengthen $p$-adic non-density to Zariski using tools from transcendence theory.
In \cite[Section 9]{LV18} it is shown that in certain moduli spaces of hypersurfaces, the integral points are not Zariski dense.
The key input is a recent transcendence result for period mappings, due to Bakker and Tsimerman \cite{BT18}.
This opens the possibility that one might prove finiteness of integral points by an inductive approach:
taking $X'$ to be the Zariski closure of the integral points in $X$, 
one would hope to use the method of \cite{LV18} to prove that, if $\dim X' \geq 1$, the integral points cannot be dense in $X'$.
To make this work one would need uniform control on the monodromy of the given family, restricted to all subvarieties $X' \subseteq X$.

\section{The method of Lawrence--Venkatesh: Finiteness}
\label{sec:LV1}

In this section, we will discuss the main ideas of the approach towards Mordell's conjecture due to Lawrence and Venkatesh. For simplicity, our main focus will be to explain the method in the case of $X = \PP^1 \backslash \{ 0,1,\infty \}$ where the proof is especially simple. Finally, we make some comments about the obstacles one faces in making this approach effective, in the example of the $2$-unit equation. 

\par Recall from the introduction that we start by constructing a map $\rho$ which attaches a Galois representation to any point on $X$. In the method of Lawrence--Venkatesh, the map $\rho$ arises from the cohomology of the fibres of a certain Parshin family $\cC \lra X$, see \S\ref{subsec:LV-General}. In the case of $X = \PP^1 \backslash \{ 0,1,\infty \}$, which we discuss first, this family is a simple modification of the classical \textit{Legendre family} of elliptic curves. 

\subsection{The $S$-unit equation. }To explain some of the ideas in the proof, we discuss the case of the $S$-unit equation in more detail. This has the benefit of being substantially simpler, while still containing many of the main ideas that go into the proof of the Mordell conjecture. To illustrate the ideas of the proof, we will start with a version of the Parshin family for which the period map $\per$ fails to be finite-to-one. Then we will give a correct argument, in which a nontrivial Galois action on $\HH^0$ of the fibers supplies the key missing ingredient.

\par Take $K = \Q$ and $S$ a finite set of primes. We denote the set of $S$-units by $\cO_S^{\times}$ and will consider the \textit{$S$-unit equation} given by
\begin{equation}
x + y = 1, \qquad x,y \in \cO_S^{\times},
\end{equation}
whose solution set is finite by Siegel's theorem. This statement represents an attractive
toy case for the Mordell conjecture; its geometric proof along the lines sketched above takes place on $X = \PP^1 \backslash \{0,1,\infty \}$. Note that we may enlarge $S$ without loss of generality, so that we may as well assume that $S$ contains $2$. 

\par The role of the Parshin family is played by the \textit{Legendre family} over $\mathcal{O}_S$. Denoting $x$ for the coordinate on $\cX = \PP^1_{\mathcal{O}_S} \backslash \{0,1,\infty \}$, this family $\cC \lra \cX$ is given by the equation
\begin{equation}
\cC \ : \ w^2 = z(z-1)(z-x).
\end{equation}
This family gives us a Galois representation $\rho(x)$ on the \'etale cohomology group $\HH^1_{\et}(\cC_{\overline{x}}, \Q_p)$, where $p$ is a prime not below any places in $S$, which is unramified in $K$. This gives the following diagram:
\begin{equation}
\label{eqn:diagram-Sunit}
\begin{tikzpicture}[->,>=stealth',baseline=(current  bounding  box.center)]
 \node[] (X) {$\cX(\cO_S)$};
 \node[right of=X, node distance=3.5cm]  (Xp) {$\cX(\cO_v)$};   
  \node[below of=X, node distance=1.5cm]  (Hf) {$\{ \mbox{Iso classes} \ \rho \}$}; 
 \node[right of=Hf,node distance=3.5cm] (Hfp) {$\{ \mbox{Iso classes} \ \rho_v \}$};
 \node[right of=Hfp,node distance=3cm](Dieu) {$\mathrm{MF}^{\phi} / \simeq$ \, .};

 \path (X)  edge node[left]{\footnotesize $\rho$} (Hf);
 \path (Xp) edge node[left]{\footnotesize $\rho$} (Hfp);
 \path (X)  edge (Xp);
 \path (Hf) edge node[above]{ } (Hfp);
 \path (Hfp) edge node[above]{\footnotesize $\D_{\cris}$} (Dieu);
 \path (Xp)  edge[dashed, bend left=20] node[above right]{\footnotesize $\per$} (Dieu);
\end{tikzpicture}
\end{equation}
Let us now make a first attempt to deduce finiteness from the above diagram. There are two major considerations to the strategy, corresponding to \textit{global} and \textit{local} aspects. The local considerations revolve around a careful analysis of the period map, via a monodromy calculation. 

\par \textbf{a. Global representations. }The Mordell conjecture will ultimately be reduced to a finiteness statement about a certain set of global Galois representations, due to Faltings. More precisely, the proof of \cite[Satz 5]{Fal83} deduces the following consequence from the classical theorem of Hermite--Minkowski:
\begin{lemma}
\label{lemma:finiteness}
Fix integers $w,d \geqslant 0$, and fix a number field $K$ and a finite set $S$ of primes of $\cO_K$.  
There are, up to conjugation, only finitely many semisimple Galois representations $\rho: G_{K} \rightarrow
\GL_d(\Q_p)$ such that
\begin{itemize}
\item[(a)] $\rho$ is unramified outside $S$, and
\item[(b)] $\rho$ is pure of weight $w$, i.e.
  for every prime $\mathfrak{p} \notin S$ all the eigenvalues of Frobenius at $\mathfrak{p}$
are algebraic integers, all of whose conjugates have complex absolute value $|\cO_K/\mathfrak{p}|^{w/2}$. 
\end{itemize} 
\end{lemma}
It should be noted that this does not make the approach of Lawrence--Venkatesh depend on the work of Faltings in an essential way, as this lemma is comparatively simple in Faltings' overall argument.

\par The semisimplicity hypothesis in Faltings' lemma is essential: there can be infinitely many nontrivial extensions between Galois representations. \footnote{As we will see in the next section, the existence of families of non-trivial extensions of a fixed set of Galois representations is precisely what underlies the method of Chabauty--Kim.} In fact, Faltings shows that all the representations we consider---which arise as subquotients of the \'etale cohomology of a curve---are semisimple. This fact requires the full weight of Faltings' argument in \cite{Fal83}. In order to give an independent proof of Mordell's conjecture, it is necessary to contemplate the possibility that some of these representations might not be semisimple. In potential algorithmic applications, we know this situation cannot arise, so we will content ourselves here with mentioning that in \cite{LV18} this is addressed by showing that all but finitely many representations in our family must be simple. This is a consequence of results of the following form:
\begin{enumerate}
\item If the global representation $\rho(x)$ has a (global) subrepresentation, then the local representation must be of a particularly special form. 
\item There are finitely many $x$ in $X(\Q_p)$ where the local representation $\rho(x)$ takes this special form.
\end{enumerate}

\par \textbf{b. The period map. }The more subtle points of the argument of Lawrence--Venkatesh lie in the study of the period map $\per$, where one systematically enlarges the base field to gain control over the Frobenius centralizers. Let us explain the need for this step, by first approaching the problem naively using the unadjusted Legendre family above. 

\par Recall that we want to show finiteness of the set of solutions to the $S$-unit equation. Since we already established the finiteness of the set of isomorphism classes of global representations $\rho(x)$ that can arise, it is tempting to try and show that the fibres of the period map $\per$ are finite. However, this is \textbf{not} true: The filtered $\phi$-modules that arise in the image of $\per$ necessarily are of the form $\HH^1_{\dR}(\cC_x, \Q_p)$, and on every good residue disk of $X(\Q_p)$ the Frobenius operator $\phi$ has a constant characteristic polynomial
\begin{equation}
f = aT^2 + bT + c \quad \in \Z_p[T]
\end{equation}
which has two roots in $\C_p$ whose valuations sum up to $1$. The number of residue disks is finite, and for each of these finitely many polynomials $f$, the filtered $\phi$-module belongs to a finite number of possible isomorphism classes, which is most easily seen with a simple case-by-case analysis:
\begin{itemize}
\item If $f$ is irreducible, then we may pick a basis $e_1$ for $\mathrm{Fil}^1$ and set $e_2 = \phi(e_1)$.  Then $\{e_1,e_2\}$ is a basis for $\HH^1_{\dR}$. With respect to this basis, we have $\mathrm{Fil}^1 = \langle e_1 \rangle$ and 
\[
\phi = \left(\begin{matrix}
0 & -b \\
1 & -a 
\end{matrix}\right).
\]
\item If $f$ is reducible, then it must have distinct roots of valuations $0$ and $1$, corresponding to eigenvectors $e_1,e_2$ which necessarily span $\HH^1_{\dR}$. Then we either have $\mathrm{Fil}^1 = \langle e_1 \rangle$ or $\langle e_2 \rangle$, or we can rescale the eigenvectors to obtain $\mathrm{Fil}^1 = \langle e_1 + e_2 \rangle$. 
\end{itemize}
In conclusion, we see that there is only a \textit{finite number} of possible isomorphism classes of filtered $\phi$-modules attached to the representations $\rho(x)$, and therefore the period map appearing in \eqref{eqn:diagram-Sunit} cannot possibly have finite fibres! Furthermore, we see from this discussion exactly what the problem is, since we had in each case so much freedom in choosing our basis, so as to move around the Hodge filtration $\mathrm{Fil}^1$ while respecting the Frobenius operator. 

\par We can rephrase the problem as follows. Fix a pair $(V, \phi)$ of a two-dimensional vector space and linear endomorphism; in our situation, $(V, \phi)$ will arise as the crystalline cohomology $\HH^1_{\cris}(\cC_x / \Z_p)$, which only depends on the reduction of $x$ modulo $p$. The possible filtrations $\mathrm{Fil}^1 \subseteq V$ are classified by the Grassmannian $\mathrm{Gr}(\mathrm{Fil}^1 \subseteq V)$. The centralizer $Z(\phi)$ acts on $\mathrm{Gr}(\mathrm{Fil}^1 \subseteq V)$, and the orbits of this action are in bijection with isomorphism classes $(V, \phi, \mathrm{Fil}^1)$ of filtered $\phi$-module with underlying $\phi$-module $(V, \phi)$. In the setting just described, $Z(\phi)$ has a Zariski-dense orbit on $\mathrm{Gr}(\mathrm{Fil}^1 \subseteq V)$, so most such filtered $\phi$-modules belong to a single isomorphism class.

\par \textbf{Interlude: Semilinearity. }
Let us take a short break to recall some crystalline theory. So far we have been applying $p$-adic Hodge theory, in particular the crystalline comparison theorem, to schemes $\cC_x$ over $\Q_p$.
In general, suppose $L_p$ is an unramified extension of $\Q_p$, and $\cC_x$ is a scheme over $L_p$,
admitting a smooth model over $\cO_{L_p}$. Then $L_p$ is Galois over $\Q_p$, with cyclic Galois group generated by a Frobenius element $\operatorname{Fr}$ that acts as the $p$-th power map on the residue field. The crystalline-de Rham cohomology $\HH^1_{\dR}(\cC_x / L_p)$ has the structure of a filtered $\phi$-module, where $\phi$ is now a \emph{semilinear} operator: it satisfies
\begin{equation}
\phi(\lambda v) = \operatorname{Fr}(\lambda) \phi(v). 
\end{equation}

\par This is important because semilinear automorphisms have small centralizers:
it's not easy for an automorphism of $V$ to both respect the action of $L_p$ and commute with $\phi$. This is made precise in the following lemma, which was proved in Lawrence--Venkatesh \cite[Lemma 2.1]{LV18}.

\begin{lemma} \label{semilinear}
Let $L_p$ be an unramified extension of $\Q_p$ of degree $e$, and let $\operatorname{Fr}: L_p \rightarrow L_p$ be the Frobenius endomorphism that acts as the $p$-th power map on the residue field. 
Let $V$ be an $L_p$-vector space of dimension $d$, and $\phi: V \rightarrow V$
a $\operatorname{Fr}$-semilinear automorphism.  Define the centralizer $\mathrm{Z}(\phi)$ of $\phi$ 
in the ring of $L_p$-linear 
endomorphisms of $V$ via
\[
  \mathrm{Z}(\phi) = \{ f: V \ \lra \ V \mbox{ an $L_p$-linear map}, \ \ f \phi = \phi f\};
\]
it is a $\Q_p$-vector space. Then 
\[
\dim_{\Q_p} \mathrm{Z}(\phi) = \dim_{L_p} \mathrm{Z}(\phi^e),
\]
where $\phi^e: V \rightarrow V$ is now $L_p$-linear. In particular, $\dim_{\Q_p} \mathrm{Z}(\phi) \leqslant (\dim_{L_p} V)^2$. 
\end{lemma}

\par \textbf{c. Finiteness. }Armed with this tool, we now return to the failed finiteness argument above, and take advantage of semilinearity to resolve the issues we were having. More precisely, we bound the size of Frobenius centralizers by considering instead the modified Parshin family
\begin{equation}
E  : w^2 = z (z-1) (z-t), \qquad t^8 = x. 
\end{equation}
For every field $K$ and $x$ in $X(K)$, the fiber $E_x$ is a geometrically disconnected curve whose $\HH^0$ is the algebra $K[t] / (t^8 - x)$. Suppose $K = \Q_p$ and $x$ is a unit in $\Q_p$ which is not a square\footnote{It is enough to consider $x$ of this form by an elementary argument based on the fact that, if $x$ is both a square and a solution to the $S$-unit equation in some number field $K$, then $\pm \sqrt{x}$ satisfy the $S$-unit equation as well. However, this does necessitate some care in the choice of $p$.}. Then $E_x$ is a curve defined over $L_p = \Q_p[x^{1/8}]$, the degree-8 unramified extension of $\Q_p$. We want to show that the map 
\[
\per : \cX(\cO_v) \ \lra \ (\mathrm{MF}^{\phi} / \simeq)
\]
is finite-to-one. On each $p$-adic residue disk $\Omega_v \subseteq \cX(\cO_v)$, the $\phi$-module $(V, \phi) = \HH^1_{\cris}(E_x)$ is constant; only the filtration varies. Thus we can regard $\per$ as a map 
\[
\per : \Omega_v \ \lra \ \mathrm{Gr}(\mathrm{Fil}^1 \subseteq V) \ \lra \ \mathrm{Gr}(\mathrm{Fil}^1 \subseteq V) / Z(\phi) \,.
\]
Since $\per$ is an analytic map from a one-dimensional source, to show it is finite-to-one, we need only show that it is not constant; in other words, that the image of $\Phi_p \colon \Omega_v \rightarrow \mathrm{Gr}(\mathrm{Fil}^1 \subseteq V)$ is not contained in a single orbit of $Z(\phi)$. This follows from the following two results:
\begin{enumerate}
\item Every orbit of $Z(\phi)$ has positive codimension in $\mathrm{Gr}(\mathrm{Fil}^1 \subseteq V)$.
\item The image of $\Phi_p$ is Zariski dense.
\end{enumerate}
The first of these two follows from the bound in Lemma \ref{semilinear}; the second, from a complex monodromy calculation. It is essential that $L_p$ have large degree over $\mathbf{Q}_p$, which comes from the assumption that $x$ is not a square in $\mathbf{Q}_p$.

\par The Zariski density of the image of $\Phi_p$ is obtained by comparing it with the complex period map $\Phi_{\C}$. Let's recall the construction of $\Phi_{\C}$.  The family $E$ of elliptic curves over $X$ gives rise to a variation of Hodge structure on $X$. Let $\Omega_{\C}$ be a contractible open subset of $X^{\an}$, containing some basepoint $x_0$ in $X(K)$, for $K$ a number field. Over $\Omega_{\C}$, the family $E$ splits as the disjoint union of eight families of elliptic curves $E^{(1)}, \ldots, E^{(8)}$. (The monodromy action of $\pi_1(X)$ preserves the splitting but permutes the eight components.) Choose an integral basis $\mathbf{B}$ for the fiberwise Betti cohomology of each elliptic curve $V^{(i)}_{\C} = \HH^1_B(E_{x_0}^{(i)})$ over $x_0$. With respect to this basis, the Hodge filtration is described by a map
\[ 
\Phi_{\C} \colon \Omega_{\C} \ \lra \  \prod_{i=1}^{8} \mathrm{Gr}(\mathrm{Fil}^1 \subseteq V^{(i)}), \]
where the Grassmannian classifies one-dimensional subspaces of the two-dimensional $V^{(i)}$.

\par The importance of $\Phi_{\C}$ to us comes from the fact that $\Phi_{\C}$ and $\Phi_p$ are, in a suitable sense, the same. (See \cite[Section 3.4]{LV18} for details.) Both period maps satisfy the same algebraic differential equation, coming from the Gauss--Manin connection. It follows that in suitable local coordinates, the (complex) power series representation of $\Phi_{\C}$ and the ($p$-adic) power series representation of $\Phi_p$ both have all their coefficients in the number field $K$, and the two power series agree. This means we can compare the images of the two period maps, and Lemmas 3.1 and 3.2 of \cite{LV18} yield the following result:
\begin{lemma}
The image of $\Phi_{\C}$ is Zariski dense if and only if the image of $\Phi_p$ is Zariski dense.
\end{lemma}

\par The advantage of this result is that establishing the Zariski-density of the map $\Phi_{\C}$ boils down to an explicit monodromy calculation, see \cite[Eqn. 3.11]{LV18}.
\begin{lemma} 
If the image of the monodromy representation of $E$ contains a Zariski-dense subset of $\operatorname{Sp}_2^d$, then the image of $\Phi_{\C}$ is Zariski dense in $\mathrm{Gr}(\mathrm{Fil}^1 \subseteq V)$.
\begin{proof}
Let $\widetilde{X}$ be the universal cover of $X$, and extend $\Phi_{\C}$ to a map
\[ 
\Phi_{\C} \colon \widetilde{X} \ \lra \ \mathrm{Gr}(\mathrm{Fil}^1 \subseteq V). 
\]
This map $\Phi_{\C}$ is $\pi_1(X)$-equivariant, where $\pi_1(X)$ acts on the Grassmannian through the monodromy representation. Since the image of monodromy is Zariski dense, the extended $\Phi_{\C}$ has Zariski-dense image. By analytic continuation, the restriction of $\Phi_{\C}$ to $\Omega_{\C}$ has Zariski-dense image as well.
\end{proof}
\end{lemma}

\begin{lemma}
\label{eqn:big-monodromy-Sunit}
The image of the monodromy representation
\[ 
\pi_1(X, x_0) \ \lra \ \operatorname{Aut} \left( \, \prod_{i=1}^8\,  \HH^1_B(E_{x_0}^{(i)}) \right) \]
contains a Zariski-dense subset of $\operatorname{Sp}_2(\Z)^8$.
\begin{proof}
This is a calculation in classical topology, see \cite[Lemma 4.3]{LV18}.
\end{proof}
\end{lemma}

\subsection{The Mordell conjecture over general $K$.}  
\label{subsec:LV-General}
After our discussion of the $S$-unit equation, we now make a brief foray into the general case, and outline how to adapt this argument to prove Mordell's conjecture. Suppose $X$ is a smooth projective curve of genus at least $2$ over $K$. We will define the \emph{Parshin family} over $X$, implicitly dependent on a parameter $q$. It will replace the Legendre family in the $S$-unit argument.

\par Let $q \geq 3$ be a prime number, and let $\Affq$ be the non-abelian group of affine-linear transformations $x \mapsto ax + b$ over $\F_q$. The action of $\Affq$ on $\F_q$ realizes $\Affq$ as a subgroup of the symmetric group $S_q$. Note also that $\Affq$ surjects onto $\F_q^{\times}$.
\begin{definition}
Let $X$ be a curve over $K$, and $x \in X(K)$ a point of $X$. An $\Affq$-cover of $X$, branched at $x$, is a curve $Z$ and a map $Z \rightarrow X$, satisfying the following properties.
\begin{itemize}
\item $Z \rightarrow X$ is \'etale over $X - \{x\}$, but not \'etale over $x$.
\item $Z \rightarrow X$ is of degree $q$.
\item For any choice of basepoint $x_0$, and for an appropriate identification of the fiber over $x_0$ with $\F_q$, the monodromy map $\pi_1(X, x_0) \rightarrow S_q$ corresponding to the cover $Z$ has image $\Affq$.
\end{itemize}
\end{definition}
For every $x$ in $X(K)$, there are finitely many isomorphism classes of $\Affq$-covers $Z \rightarrow X$ branched at $x$. The Parshin family $Y \rightarrow X$ is characterized by the property that the fiber $Y_x$ is geometrically the disjoint union of these finitely many curves.

\par In the $S$-unit argument, the key semilinearity bound came from taking an $8$-th root of $x$
(along with the elementary assumption that $x$ is not a square). Here the corresponding bound comes from the torsion on the Jacobian $J$ of $X$, which is guaranteed to have a nontrivial Galois structure. Specifically, for any $\Affq$-cover of $X$, the composed map
\[ 
\pi_1(X - \{x\}) \ \lra \ \Affq \ \lra \ \F_q^{\times}
\]
gives a degree-$(q-1)$ cover of $X$ that turns out to be unramified everywhere, even over $x$.
This cover in turn corresponds to a $(q-1)$-torsion point on the dual of the Jacobian. We choose $q$ and $p$ so that the Frobenius at $p$ acts with sufficiently large orbits on $J[q-1]$; this in turn guarantees that the components of each fiber $Y_x$ are defined over large $p$-adic fields, so we can leverage the semilinearity lemma \ref{semilinear}.

\par As with the $S$-unit equation, a calculation in the classical topology is needed to show that the Parshin family has big monodromy.
Fix $X$ and $x$,  and let $Z_1, \ldots, Z_N$ be the $\Affq$-covers of $X$ branched at $x$.
We want to determine the image of the monodromy action
\[ 
\operatorname{Mon} \colon \pi_1(X, x) \ \lra \ \operatorname{Aut} 
\left( \, \prod_{i} \, \HH^1_B(Z_i) \right)
\]
as an algebraic group. The cohomology of each $Z_i$ contains a copy of $\HH^1_B(X)$; define\footnote{The symbol ``Pr'' stands for ``primitive.''} 
\[
\HH^1_{\rm Pr}(Z_i) = \HH^1_B(Z_i) / \HH^1_B(X).
\]
The map $\operatorname{Mon}$ descends to an automorphism of $\prod_i \HH^1_{\rm Pr} (Z_i)$. We need the following big monodromy result.
\begin{theorem}
\label{big_mon}
The Zariski closure of the image of
\[  
\operatorname{Mon} \colon \pi_1(X, x) 
\ \lra \ 
\operatorname{Aut} \left( \prod_{i}  \HH^1_{\rm Pr}(Z_i) \right) 
\]
contains the group
\[ 
\prod_{i} \operatorname{Sp} ( \HH^1_{\rm Pr}(Z_i)).
\]

\end{theorem}
This theorem is really saying that the image of monodromy is as big as possible: we know for abstract reasons that the identity component of the Zariski closure of the image is no larger than the product of symplectic groups. We say a few words about the main ideas that go into the proof:

\par The monodromy action of $\pi_1(X, x)$ on the covers $Z_i$ extends to an action of the full mapping class group\footnote{The mapping class group is the group of \emph{topological} automorphisms of the topological surface $X$ fixing the point $x$, up to isotopy fixing $x$.  The book of Farb and Margalit \cite{FM12} is an excellent introduction and reference on mapping class groups.} $\MCG(X - \{x\})$. By the Birman exact sequence, $\pi_1(X, x)$ is a normal subgroup of $\MCG(X - \{x\})$. Since the symplectic group is simple modulo center, we can deduce Theorem \ref{big_mon} if we know that the Zariski closure of the image of
\[  
\operatorname{Mon} \colon \MCG(X - \{x\}) \ \lra \ \operatorname{Aut} \left( \prod_i \HH^1_{\rm Pr}(Z_i)\right) 
\]
contains said product of symplectic groups.
The benefit to working with the full mapping class group is that we now have access to Dehn twists,
a particularly simple class of automorphism that is amenable to explicit calculation. Dehn twists map to unipotent automorphisms via $\operatorname{Mon}$, and the proof concludes by producing a collection of unipotent automorphisms that generates the full symplectic group.

\par The study of mapping class group representations like $\operatorname{Mon}$ is a big subject in geometric topology. Looijenga \cite{Looijenga} studied the analogous question for abelian covers.
Grunewald, Larsen, Lubotzky, and Malestein \cite{GLLM} study (unramified) covers of
\emph{compact} surfaces, and in a recent paper~\cite{ST} Salter and Tshishiku study covers whose covering group is the Heisenberg group. These results are stronger than ours: they all show that the image of the representation has finite index in an appropriate arithmetic group, rather than merely being Zariski dense.

\section{The method of Lawrence--Venkatesh: Effectivity}
\label{sec:LV2}

We now discuss the extent to which we expect the work of Lawrence--Venkatesh to yield a method for \textit{explicitly determining} the set $X(K)$ in examples. Since it is so recent, it is unsurprising that this aspect of the method of Lawrence--Venkatesh does not yet seem to be addressed in the literature. In this section we adopt a more speculative tone, merely making some brief comments about various ingredients that would likely be needed to parlay this method into an algorithm for bounding the number of rational points on a curve over a number field; which would yield, in a weak sense, a form of ``algorithmic Mordell.''

\par Roughly speaking, a potential form of such a hoped-for algorithm is as follows.

\begin{algorithm}
Take as input a number field $K$, a smooth projective curve $X$ over $K$, and a power $v^n$ of a good\footnote{The method of LV requires $p$ to satisfy a certain Galois-theoretic condition; here we will simply call primes satisfying that condition ``good'' primes. The condition is needed to guarantee that a certain extension of $K_p$ is of large degree, and is analogous to the requirement in the $S$-unit equation above that $x$ not be a square in $\mathbf{Q}_p$. Choosing a good $p$ presents no algorithmic difficulty.} prime ideal $v$ of $\mathcal{O}_K$. Return as output a finite list of points\footnote{possibly with multiplicities} in $\cX(\cO_v)$, to any desired finite precision which is guaranteed to include all the rational points of $X$.
\end{algorithm}

\par It should be mentioned that such an algorithm, until an efficient implementation proves the contrary, is at risk of being prohibitively slow so as to be useless from a practical standpoint. The essential difficulty lies in enumerating Galois representations with prescribed ramification; modularity results for the representations in question, if known, could speed up the algorithms significantly. One possible approach to the calculation is proposed in what follows. It has four essential components, each of which we briefly discuss below. It should be noted that whereas many of the separate ingredients have been extensively studied in the literature, the method of Lawrence--Venkatesh has so far not been made effective, and therefore the ideas in this section are tentative. It would be very interesting to explore the effectivity of this method further, and make a serious attempt at a computational version of this method. 

\begin{remark}An algorithm of the above form may return extraneous points, not corresponding to a rational point. This phenomenon arises also in Chabauty's method, though in the example in \S\ref{sec:Intro}  it was circumvented by exhibiting two independent analytic sets, which was possible since $g-r = 2$. Likewise, it is conceivable that one can circumvent in the method of Lawrence--Venkatesh by varying the choice of $q$ in the covering group $\Affq$. Alternatively, one could attempt to apply the Mordell--Weil sieve, see \S\ref{subsec:MordellWeil}.\end{remark}

\subsection{Enumerating global Galois representations. }

Faltings's finiteness lemma for Galois representations (Lemma \ref{lemma:finiteness}) can be made effective; we expect this to be the most computation-intensive part of the algorithm. Recall that we want to enumerate all global Galois representations
\[ \rho: G_{K} \ \lra \ \GL_d(\Z_p) \]
that could arise from our family, in the sense of Lemma \ref{lemma:finiteness}. We know the following about $\rho$:
\begin{itemize}
\item We are given a finite set $S$ of places of $K$, outside of which $\rho$ is unramified.
\item For every prime $\mathfrak{p} \notin S$, all the eigenvalues of Frobenius at $\mathfrak{p}$ are Weil numbers of weight 1/2.
\item The representation $\rho$ is semisimple.
\end{itemize}

On the one hand, we can list all possible mod-$p^n$ representations for any $n$. First, one enumerates all possible residual representations
\[ 
\rho_1: G_{K} \ \lra \ \GL_d(\F_p). 
\]
This is a straightforward application of Hermite--Minkowski finiteness. The residual representation has finite image, so it corresponds to an extension $L_1$ of $K$ of degree at most $\left | \GL_d(\F_p) \right |$. The ramification condition translates to a bound on the discriminant of $L_1$. One can find all possible number fields $L_1$ by a targeted Hunter search \cite[\S 9.3]{Coh00}.
However, the time complexity of such a search (for fixed $K$ and $S$) is doubly exponential in the degree $[L_1: K]$, so it may be necessary to further refine the search using more specifics of the situation at hand.

\par Second, for each residual representation, one lifts successively to mod-$p^n$ representations
\[ 
\rho_n: G_{K} \ \lra \ \GL_d(\Z/p^n),  
\]
which correspond to a tower of fields $L_n$.
The successive extensions $L_{n+1} / L_n$ are abelian, so they can be found by class field theory.
(Everything we need from class field theory can be done algorithmically; see \cite{Coh00}.)
To do this, we need to compute ideal class groups and unit groups of number fields whose degrees grow exponentially in $n$; this is again a computationally expensive task.

On the other hand, given the residual representation $\rho_1$,
the Faltings--Serre method (see for example \cite{Deligne_ladic})
allows one to compute effectively a finite set of primes $\mathfrak{p}_1, \ldots, \mathfrak{p}_s$ such that
for any semisimple $\rho$ lifting $\rho_0$, the \emph{rational} representation
\[
G_{K} \ \lra \ \GL_d(\Q_p)
\]
 is determined by the Frobenius traces 
\[ 
\Tr ( \operatorname{Fr}_{\mathfrak{p}_i} | \rho) 
\]
at these finitely many primes. (In general, there may be multiple isomorphism classes of \emph{integral} representation,  as the rational representation may have more than one stable $\Z_p$-lattice.) The condition on Frobenius eigenvalues guarantees that there are only finitely many possible values for $ \Tr ( \operatorname{Fr}_{\mathfrak{p}_i} | \rho)$, for each $i$. We can choose $n_0$ large enough that, for each fixed $i \in \{1, \ldots, s\}$, no two of these possible values are congruent modulo $p^{n_0}$. Then any mod-$p^{n_0}$ representation $\rho_{n_0}$ can lift to at most one semisimple $p$-adic representation.

\par The strategy, then, is as follows.  First, make a list of all (finitely many) possible tuples
\[ 
(\Tr ( \operatorname{Fr}_{\mathfrak{p}_i} | \rho))_{i \in \{1, \ldots, s\}}; 
\]
we'll call such a tuple a \emph{candidate}. As described above, we can enumerate all mod-$p^n$ representations for some $n \geq n_0$. We compute their Frobenius traces and match them with candidates, discarding candidates that don't match any representation, and vice-versa. We can repeat this procedure for any desired $n$; the list of candidates will get shorter, as spurious candidates are deleted.

\subsection{Computing the Parshin family. }Before we get to the purely local part of the computation, which consists of describing the $p$-adic period map $\per$, we are faced with the problem of finding an explicit set of algebraic equations defining the Parshin family 
\[
\cC \ \lra \ X,
\]
whose fibres are finite covers of $X$ branched over the variable point $x$. This is an instance of the Riemann--Hurwitz problem. Computational work on branched covers of curves is particularly well-developed in the case of Bely\u{\i} covers of $\PP^1$; see \cite{SV14} for an overview. The covers appearing in our setting are solvable, and we expect that explicit calculations on the Jacobian could provide a fruitful approach. 

\par The solvability of the covering group $\Affq$ has the following geometric
interpretation.  Suppose $Z \rightarrow X$ is an $\Affq$-cover, branched at $x$.  Let
$Z^{\Gal}$ be the Galois closure of $Z$; this is a cover of $X$ of degree $q(q-1)$,
ramified only above $x$ and having Galois group $\Affq$.  The quotient map $\Affq
\rightarrow \F_q^\times$ gives a curve $Z^{\Ab}$, corresponding by the Galois
correspondence to $\F_q^\times$.  Thus we have the tower of covers
\[ Z^{\Gal} \ \lra \ Z^{\Ab} \ \stackrel{\pi}{\lra}\  X. \]
In this tower, $Z^{\Ab}$ is an unramified abelian cover of $X$ of degree $q-1$, and $Z^{\Gal}$ is an abelian cover of $Z^{\Ab}$ of degree $q$, ramified at exactly the points of $\pi^{-1}(x)$. The curve $Z$ can be recovered as a quotient of $Z^{\Gal}$.

\par This suggests the following strategy to compute the Parshin family $Y$, each of whose fibers is a union of $\Affq$-covers $Z$. First, we attempt to compute abelian covers (both unbranched and branched) of arbitrary curves, by finding torsion points on algebraic generalized Jacobians. To describe one strategy\footnote{An alternative approach to computing covers of curves is by Hensel lifting from a finite field, as in \cite{Mascot}. }, we will restrict attention to unramified covers and the (ungeneralized) Jacobian. In this setting, we want to find a divisor $D$ on the curve $X$, along with a meromorphic function $f$ on $X$ such that
\[ 
\mathrm{div}(f) = r D,
\]
which amounts to looking for $r$-torsion on the Jacobian of $X$. The Jacobian has an algebraic incarnation as a variety classifying divisor classes on $X$ and an analytic incarnation as a complex torus. It is of course trivial to identify torsion points on the \emph{analytic} Jacobian;
what we need is to describe them as points on the algebraic Jacobian.

\par Fix a basepoint $b \in X(\mathbf{C})$. By integration we can compute coordinates on the analytic torus $\operatorname{Jac} X$, along with the analytic Abel--Jacobi map
\[ 
X \ \lra \ \operatorname{Jac} X. 
\]
In the other direction, let $g$ be the genus of $X$.  We want to invert the map
\begin{equation}
\label{eqn:symX}
\operatorname{Sym}^g X \ \lra \ \operatorname{Jac} X, 
\end{equation}
to realize a point of the analytic Jacobian as a divisor on $X$. This map is a birational equivalence, but not an isomorphism. On a Zariski-dense subset of $\operatorname{Jac} X$, the map can be inverted, for example, by theta function methods \cite[Theorem II.3.1]{Tata1}, by Puiseux series methods \cite{CMSV}, or by computations in Grassmannians arising from Riemann-Roch theory \cite{KM03, CMSV}. 
A general algorithm appears in \cite[\S 3.3]{CMSV}.

\par If we can compute arbitrary abelian covers, we could try to determine all the covers $Z^{\Gal}$ for any fixed point $x$; from there one computes $Z$ by Galois theory on the function field. In other words, we can compute the fiber $Y_x$ of the Parshin family over any given point $x \in X$. To compute the Parshin family as an algebraic family, we are faced with the need to interpolate these fibers, perhaps by Puiseux series methods.

\subsection{Computing the $p$-adic period map. }We now come to the local part of the computation, where a description of the $p$-adic period map $\per$ reduces to a computation with $p$-adic cohomology in families. There is a vast literature on this subject, and this step is therefore likely to be more accessible and efficient than the others\footnote{Indeed, the algorithms mentioned here are crucial ingredients for the effective method of Chabauty--Kim, as we will see in \S\ref{sec:Examples}. }. We give a brief overview of some results in the literature, for more detailed treatments that address also the history of the subject, see Kedlaya \cite{Ked05,Ked07}.

\par The basic problem is the following: Suppose we are given a curve $\cC_x$ over a $p$-adic field $K_v$ and want to compute the filtered $\phi$-module structure of $\HH^1_{\dR}(\cC_x/K_v)$. Representing this space by differentials of the second kind, the Hodge filtration is easily worked out, and it is the Frobenius operator $\phi$ that forms the essence of this problem. When $\cC_x$ is hyperelliptic, Kedlaya \cite{Ked01} introduced an efficient algorithm, a variant of which we will see in action for the examples of the genus $2$ curves  in \S\ref{sec:Examples}. There are two main ingredients for the computation: 
\begin{itemize}
\item An appropriate lift of Frobenius on the functions in a ($p$-adic analytic) open subset of $\cC_x$,
\item A reduction algorithm in de Rham cohomology, that writes an arbitrary differential as the sum of an exact differential and a linear combination of our basis differentials. 
\end{itemize}
By applying the reduction in cohomology to the image of a set of basis differentials under this Frobenius lift, we may obtain a matrix of the Frobenius operator $\phi$, up to some precision $v^n$.

\par This method has seen extensive developments since \cite{Ked01}, notably by Lauder \cite{Lau04,Lau06} who introduced the \textit{fibration method}. This method makes use of the Frobenius structure on the sheaf of relative $q$-th de Rham cohomology $\cH^{q}_{\dR}(X/S)$ of a smooth morphism $X\lra S$ between smooth varieties over $K_v$. The variation of the de Rham cohomology of the fibres in this family is described by the Gau{\ss}--Manin connection
\[
\nabla_{\GM} : \cH^q_{\dR}(X/S) \ \lra \Omega^1_{X/S} \otimes \cH^q_{\dR}(X/S),
\]
which gives a system of differential equations known as the \textit{Picard--Fuchs equations}, whose study was taken up in the 19\textsuperscript{th} century. Suppose we find a local lift of Frobenius $\phi$ on $S$, then the pullback of the relative de Rham cohomology $\cH^q_{\dR}(X/S)$ by $\phi$ is isomorphic \textit{as a vector bundle with connection} to the original one. In concrete terms, let us suppose that $S$ is a curve, then we may express this in matrix form as
\begin{equation}
\label{eqn:GM}
NFdt + \frac{\partial}{\partial t}F = \left(\frac{\partial}{\partial t}\phi(t)\right) F\phi(N)dt
\end{equation}
by choosing a local coordinate $t$ on $S$, and a basis of the relative de Rham cohomology, with respect to which we obtain a matrix $F(t)$ describing the Frobenius operator on the fibres, and $N(t)dt$ describing the Gau{\ss}--Manin connection. This equation is very useful. For instance, if $F(t)$ can be computed for a single value of $t=t_0$, then we may solve these $p$-adic differential equation using $F(t_0)$ as an initial condition. Lauder \cite{Lau04,Lau06} uses this idea to compute the Frobenius action in families. It is surprisingly versatile, applying both to individual curves with a map to $\PP^1$ as well as families of curves. It has been developed in many subsequent papers of which we mention the recent algorithms of Tuitman \cite{Tui16,Tui17}, and the references contained therein, which vastly extend the range of applicability of these ideas. 

\subsection{Compare the global Galois representations with the $p$-adic periods. }

We suggest two approaches. The first is to use $p$-adic Hodge theory, along the lines of Fontaine--Laffaille theory \cite{FL82}. We are given a mod-$p^n$ global Galois representation, presented as a polynomial whose splitting field is its kernel. We can determine the corresponding local representation at $p$, in terms of extensions of $\mathbf{Q}_p$. Fontaine and Laffaille define a functor $\underline{U}_S$ from a certain category of finite-length filtered $\phi$-modules to the category of Galois representations \cite[\S 0.6]{FL82}. One expects that Fontaine--Laffaille theory can be made algorithmic: given a mod-$p^n$ Galois representation, we should be able to determine whether it is in the image of this functor, and if so, describe the underlying filtered $\phi$-module. We can then compare these $\phi$-modules with the $\phi$-modules arising from the $p$-adic period map,
to determine a list of candidate points.

\par Our second approach avoids filtered $\phi$-modules entirely, by working directly with Galois representations. It is a consequence of Fontaine--Laffaille theory that the mod-$p^n$ local Galois representation $\rho_x$ depends only on the reduction of $x$ modulo $v^{n+1}$. Using this, we can compute explicitly all the possibilities for the local Galois representation at $p$, and match them explicitly with the list of ``candidates'' from the global Galois calculation. In other words, for each candidate $\rho$, we obtain a list of mod-$v^{n+1}$ points of $X$, the local representations at which agree modulo $p^n$ with $\rho$. For each of these mod-$v^{n+1}$ residue classes, we then use the period map $\per$ to compute a bound on the number of rational points in the class.

\section{The method of Chabauty--Kim: Finiteness}
\label{sec:CK1}
In this section, we discuss the approach to Mordell's conjecture due to Minhyong Kim. It follows the same pattern as the method of Chabauty--Coleman discussed in the introduction, and as such it depends on some geometric input, replacing the condition $r < g$ by something weaker, which may be done at the cost of replacing the $p$-adic Tate module $V$ by a more sophisticated quotient of the fundamental group. We discuss in some detail the particular case of a quotient arising from a geometric correspondence \cite{BD18,BD19,BDMTV19} using the geometric language of Edixhoven--Lido \cite{EL19}.

\subsection{Quotients of the fundamental group. }
\label{subsec:quotients}
To motivate an interest in unipotent quotients of the algebraic fundamental group for Diophantine applications, it is instructive to first recall the \textit{section conjecture} of Grothendieck \cite{Gro97}, which states that the map
\[ 
\begin{array}{cccc}
\rho \ : & X(\Q) &\longrightarrow& \HH^1 \! \left(G_{\Q},\pi^{\et}_1(\overline{X},b)\right), \\ [5pt]
 & x & \longmapsto & [ \ \pi^{\et}_1(\overline{X};b,x) \ ] \\
\end{array}
\]
which attaches to every rational point the class of the Galois representation defined by the corresponding \textit{path torsor} of the algebraic fundamental group, should be an isomorphism. In other words, every torsor of the fundamental group should necessarily arise from a rational point. This provides us with the tantalizing possibility of studying the set of such torsors in lieu of the set $X(\Q)$. Unfortunately, the cohomology set that classifies these torsors does not seem to have much structure with which we can work. 

\par On the other end of the spectrum, we already saw that the twists of the $p$-adic Tate module $V$ of the Jacobian $J$ of $X$, which is essentially the abelianization of the fundamental group, are classified by an object which is very closely related to $J$, and which therefore has a tremendous amount of structure. That said, this association only gives us enough information under the additional assumption that $r < g$. 

\par In summary, we could roughly describe the situation by saying that the association
\begin{equation}
\rho \ : \ X(\Q) \ \lra \ \HH^1 \! \left(G_{\Q},\pi^{\et}_1(\overline{X},b)\right)
\end{equation}
in the section conjecture has a target with \textit{too little} structure, whereas the association 
\begin{equation}
\rho \ : \ X(\Q) \ \lra \ \HH^1_f(G_{\Q},V)
\end{equation}
appearing in the method of Chabauty--Coleman has a target with \textit{too much} structure. The latter statement is meant in the sense that $\rho$ factors through the Jacobian, and in situations where $r \geq g$ this kills some crucial non-abelian information needed to understand $X(\Q)$. In the method of Chabauty--Kim, we allay the difficulties inherent to both settings by working with a suitable intermediate quotient, balancing the availability of structure on the sets $\HH^1$ against our ability to explicitly describe the target. 
We consider quotients of the fundamental group that are unipotent.\footnote{Strictly speaking, quotients of the $\Q_p$-unipotent \'etale fundamental group studied in Deligne \cite{Del89}.}  

\par The strategy for proving finiteness follows the same pattern as our discussion of the method of Lawrence--Venkatesh. First, one attempts to gain sufficient control over the set of global representations involved, and second, one studies the local representations via the analytic properties of an associated period map. 

\par \textbf{a. Global representations. } A general theorem of Kim (\cite[Proposition 2]{Kim05} and \cite[p. 118]{Kim09}) states that if $\U$ is a unipotent quotient satisfying certain technical assumptions which we will not state here, the set $\HH^1_f(G_K,\U)$ carries the structure of an algebraic variety, dubbed \textit{Selmer variety}, such that the localization map
\begin{equation}
\HH^1_f(G_{\Q},\U) \ \lra \ \HH^1_f(G_{\Q_p},\U)
\end{equation}
between the global and local Selmer varieties is algebraic. The algebraic nature of this map allows us to gain control over the image of the global Selmer variety, typically by showing that the global Selmer variety is of lower dimension than the local Selmer variety, so that the image cannot be Zariski dense. 

\par \textbf{b. The period map. }As was the case in the method of Lawrence--Venkatesh, the control of global representations can be turned into a proof of finiteness by controlling a $p$-adic period map. In the method of Chabauty--Kim, this means concretely that one establishes that the association 
\begin{equation}
\label{eqn:local_rho_U}
\rho \ : \ X(\Q_p) \ \lra \ \HH^1_f(G_{\Q_p}, \U),
\end{equation}
of the path torsor of $\U$ attached to a point, has an image which is Zariski dense. Typically, the quotient $\U$ is of a ``motivic'' nature, in which case the association $\rho$ in \eqref{eqn:local_rho_U} has a de Rham realisation
\begin{equation}
\per \ : \ X(\Q_p) \ \lra \ \mathrm{MF}^{\phi}
\end{equation}
which can be expressed as a linear combination of iterated Coleman integrals of differentials. A general theorem of Kim \cite[Theorem 1]{Kim09} establishes the linear independence of such iterated integrals, which often implies the Zariski density of the image of \eqref{eqn:local_rho_U} by $p$-adic Hodge theory. 

\par \textbf{c. Finiteness. }In conclusion, we are left with the following attractive strategy to study the set of rational points $X(\Q)$: Suppose that we can construct a specific finite-dimensional unipotent quotient $\U$ satisfying the technical hypotheses required for the representability of the Selmer varieties, such that furthermore
\begin{enumerate}
\item we can prove that $\dim \HH^1_f(G_{\Q},\U) < \dim \HH^1_f(G_{\Q_p},\U)$,
\item the quotient is ``motivic'', so that we have a $p$-adic period map
\begin{equation*}
\per \ : \ X(\Q_p) \ \lra \ \mathrm{MF}^{\phi}
\end{equation*}
which is a linear combination of iterated integrals of differentials on $X$, and
\item we can find a \textit{computable} condition on elements of the image of $\per$ to come from a point in $X(\Q)$.
\end{enumerate}
Once we manage to find a quotient $\U$ satisfying these conditions, we consider the diagram
\begin{equation}\label{eqn:diagram_kim}
\begin{tikzpicture}[->,>=stealth',baseline=(current  bounding  box.center)]
 \node[] (X) {$X(\Q)$};
 \node[right of=X, node distance=3.7cm]  (Xp) {$X(\Q_p)$};   
  \node[below of=X, node distance=1.5cm]  (Hf) {$\HH^1_f(G_{\Q},\U)$}; 
 \node[right of=Hf,node distance=3.7cm] (Hfp) {$\HH^1_f(G_{\Q_p},\U)$};
 \node[right of=Hfp,node distance=3.5cm](Dieu) {$\mathrm{MF}^{\phi}$.};

 \path (X)  edge node[left]{\footnotesize $\rho$} (Hf);
 \path (Xp) edge node[left]{\footnotesize $\rho$} (Hfp);
 \path (X)  edge (Xp);
 \path (Hf) edge node[above]{ } (Hfp);
 \path (Hfp) edge node[above]{\footnotesize $\mathbf{D}_{\cris}$} (Dieu);
 \path (Xp)  edge[dashed, bend left=20] node[above right]{\footnotesize $\per$} (Dieu);
\end{tikzpicture}
\end{equation}
The first two conditions on $\U$ are the active ingredients for deducing finiteness. The
first condition is the analogue of the condition $``r < g"$ appearing in the method of
Chabauty--Coleman, and allows us to control the image of the global Selmer variety. When
combined with a concrete understanding of the period map $\per$ provided by the second
condition (for instance, enough to show Zariski-density of~\eqref{eqn:local_rho_U}, see \cite[Theorem 1]{Kim09}) the above commutative diagram implies that $X(\Q_p)$ intersects the image of the global Selmer variety in a finite set of points. In particular, this shows that $X(\Q)$ is finite. 

\par Finding suitable quotients that satisfy the first two conditions is the subject of many works, and is typically done by considering quotients $\U$ arising from powers of the augmentation ideal, see for instance Kim \cite{Kim05,Kim09}, Coates--Kim \cite{CK10} and Ellenberg--Hast \cite{EH}. The third condition is relevant for the \textit{explicit determination} of $X(\Q)$ and will reappear later. 

\subsection{Geometric correspondences on $X$. }
\label{subsec:correspondences}
We now discuss one instance where such a quotient can be constructed, under the additional assumption that the Jacobian $J$ of $X$ has non-trivial N\'eron--Severi rank, following \cite{BD18,BDMTV19}. To offer a different perspective on the constructions in \textit{loc. cit.} we opt for the more geometric reformulation of this theory following the beautiful work of Edixhoven--Lido \cite{EL19}. It should be noted that in \cite{EL19} this geometric viewpoint is retained to find a method for the effective determination of $X(\Q)$, but in \S\ref{sec:CK2} we instead opt for the cohomological language of \cite{BD18,BDMTV19}. 

\par Recall that the N\'eron--Severi group of a smooth proper variety is the group of components of its Picard scheme. In the situation at hand, we have chosen a base point $b$ in $X(\Q)$, which gives us an associated Abel--Jacobi map $X \lra J$. By functoriality, we obtain the following diagram:
\begin{equation}
\label{eqn:diagram-UZ}
\begin{tikzpicture}[node distance=2.3cm, >=stealth', baseline=(current  bounding  box.center)]

	\node (A1) {$\mathrm{Pic}^0(J)$};
	\node (A2)[below of=A1, node distance=1.3cm] {$\mathrm{Pic}^0(X)$};
	\node (AC1)[right of=A1]{$\mathrm{Pic}(J)$};
	\node (AC2)[right of=A2]{$\mathrm{Pic}(X)$};
	\node (C1)[right of=AC1]{$\mathrm{NS}(J)$};
	\node (C2)[right of=AC2]{$\Z$};

	\node (001)[right of=C1, node distance=1.5cm]{$0$};
	\node (002)[right of=C2, node distance=1.5cm] {$0$.};
	\node (01)[left of=A1, node distance=1.5cm]{$1$};
	\node (02)[left of=A2, node distance=1.5cm] {$1$};
	
	\draw [->,above] (A1) to node {} (AC1);
	\draw [->,above] (A2) to node {} (AC2);
	\draw [->] (AC1) to node[above] {} (C1);
	\draw [->] (AC2) to node[above] {} (C2);
	\draw [-,double distance=2pt] (A1) to node[right] {} (A2);
	\draw [->] (C1) to node[right] {} (C2);
	\draw [->] (AC1) to node[right] {} (AC2);
	\draw [->] (C1) to node[above] {} (001);
	\draw [->] (C2) to node[above] {} (002);	
	\draw [->] (01) to node[above] {} (A1);
	\draw [->] (02) to node[above] {} (A2);
\end{tikzpicture}
\end{equation}
The N\'eron--Severi group $\mathrm{NS}(J)$ is a finitely generated group, of rank $\mathrm{rk}_{\mathrm{NS}}$ which is called the \textit{N\'eron--Severi rank} of $J$. Now suppose that we have a non-trivial class $Z$ in $\mathrm{NS}(J)$ which maps to zero in $\Z \simeq \mathrm{NS}(X)$ in the above diagram\footnote{Such a class always exists when $\mathrm{rk}_{\mathrm{NS}} > 1$, which is true for many examples of interest, including modular curves. }. Then, by the identification of $\mathrm{Pic}^0(J)$ with $\mathrm{Pic}^0(X)$ there is a unique lift of $Z$ to an element of $\mathrm{Pic}(J)$ which is trivial when restricted to $X$. In other words, $Z$ uniquely determines a (non-trivial) line bundle $\mathscr{L}_Z$ on $J$ which is trivial when restricted to $X$, and hence we obtain a lift of the Abel--Jacobi map
\begin{equation}
\begin{tikzpicture}[node distance=2cm, >=stealth', baseline=(current  bounding  box.center)]

	\node (X) {$X$};
	\node (J)[right of=X]{$J$.};
	\node (L)[above of=J, node distance=1.1cm]{$\mathscr{L}_Z$};
	
	\draw [->] (X) to node {} (J);
	\draw [->] (L) to node {} (J);
	\draw [->, dashed, bend left = 20] (X) to node[above] {} (L);
	
\end{tikzpicture}
\end{equation}
This lifting of the Abel--Jacobi map, or equivalently this trivialization of the line bundle $\mathscr{L}_Z$ restricted to $X$, is a priori uniquely determined up to multiplication by elements of $\Q^{\times}$. As explained in Edixhoven--Lido \cite{EL19}, one can determine it up to $\Z^{\times} = \{\pm 1\}$, and hence essentially uniquely, at the cost of taking a small\footnote{It suffices to take the least common multiple of the exponents of the N\'eron component groups at all primes of bad reduction.} tensor power of $\mathscr{L}_Z$ by spreading out the geometry over $\Z$ and working with the N\'eron model of $J$. In conclusion, we obtain an essentially unique lift
\begin{equation}
X \ \lra \ \mathscr{L}_Z^{\times} := \mathbf{Isom}_J(\cO,\mathscr{L}_Z).
\end{equation}

\par The scheme $\mathscr{L}_Z^{\times}$ is a $\G_m$-torsor\footnote{Since the class of its line bundle in $\mathrm{Pic}(J)$ maps to a non-zero element of $\mathrm{NS}(J)$, this $\G_m$-torsor is not a group.} over the Jacobian $J$. We define $\U$ to be the $\Q_p$-\'etale fundamental group of $\mathscr{L}_Z^{\times}$. This group is non-abelian, and may be understood geometrically as follows. One can show (see for instance Bertrand--Edixhoven \cite[\S~4]{BE19} for the arguments in the $\C$-analytic setting) that there is a co-final system of \'etale coverings 
\[
\pi_n : (\mathscr{L}_Z^{\times})^{\otimes n } \ \lra \ \mathscr{L}_Z^{\times} 
\]
obtained by composing the pullback of the map $[n]$ on the Jacobian with the $n$\textsuperscript{th} power map on fibres. The Galois group $\U_n$ of this \'etale cover is a central extension
\[
1 \ \lra \ \mu_n \ \lra \ \U_n \ \lra \ J[n] \ \lra \ 0,
\]
so that $\U$ is a Heisenberg group, and as a Galois representation it is an extension of $V$ by $\Q_p(1)$. Suppose that $x$ is a point in $X(K)$ for $K$ equal to $\Q$ or $\Q_p$. Then in analogy with \eqref{eqn:torsor-J} we obtain a torsor of $\U$ from the inverse limit of the preimages of $x$ under the maps $\pi_n$, i.e. 
\begin{equation}
\label{eqn:torsor-L}
\Q_p \otimes_{\Z_p} \left( \varprojlim_{n} \, \pi_n^{-1}(x)\right), \qquad \pi_n \ : \ (\mathscr{L}_Z^{\times})^{\otimes n}\ \lra \ \mathscr{L}_Z^{\times}.
\end{equation}
Such torsors are classified by the cohomology group $\HH^1(G_K,\U)$ and satisfy certain local conditions which we will not make explicit here. In conclusion, we obtain an association $\rho$ analogous to \eqref{eqn:AJ-torsor}:
\begin{equation}
\label{eqn:U-torsor}
\rho \ : \ X(K) \ \lra \ \mathscr{L}_Z^{\times}(K) \ \lra \ \HH^1_f(G_K, \U).
\end{equation}

\subsection{Finiteness of $X(\Q)$. }The quotient $\U$ attached to a N\'eron--Severi class $Z$ as above has dimension $2g+1$ as a $\Q_p$-vector space. More precisely, as a Galois representation, it is an extension of the form
\begin{equation}
0 \ \lra \ \Q_p (1) \ \lra \ \U \ \lra \ V \ \lra \ 0,
\end{equation}
where $\Q_p(1)$ is the one-dimensional representation given by the cyclotomic character. The simple nature of this one-dimensional graded piece is responsible for the proof that the quotient $\U$ satisfies the first condition on our wish list in \S\ref{subsec:quotients}. Indeed, this can be deduced from the statements
\begin{equation}
\begin{array}{lllll}
\HH^1_f(G_{\Q_{\phantom{p}}}, \Q_p(1)) &=& \Z^{\times} \widehat{\otimes} \Q_p &=& 0,  \\ [3pt]
\HH^1_f(G_{\Q_p}, \Q_p(1)) &=& \Z_p^{\times} \widehat{\otimes} \Q_p &=& \Q_p
\end{array}
\end{equation}
which result, via the simple argument in \cite[Lemma 3.1]{BD18}, in the statements
\begin{equation}
\begin{array}{ccl}
\dim \HH^1_f(G_{\Q_{\phantom{p}}},\U) & \leq & r,\\ [3pt]
\dim \HH^1_f(G_{\Q_p},\U) & = & g+1. 
\end{array}
\end{equation}

\par The quotient $\U$ is also \textit{motivic} in nature, as its geometric definition via the $\G_m$-torsor $\mathscr{L}_Z^{\times}$ shows. In particular, besides the Galois representation $\U$, there is also a de Rham realisation $\U^{\dR}$, which is a quotient of the de Rham fundamental group of $X$, see \cite{Kim05,Kim09} for more precise definitions. The theorem of Kim \cite[Theorem 1]{Kim09} discussed in \S\ref{subsec:quotients} then implies that the image of $X(\Q_p)$ under $\rho$ is Zariski dense. 

\par This allows us to deduce finiteness of $X(\Q)$ for certain curves $X$. Suppose that
$r=g$, so that we are just outside of the range where the method of Chabauty--Coleman
applies, and assume furthermore that the N\'eron--Severi rank $\mathrm{rk}_{\mathrm{NS}}$ of $J$ is at least $2$, so that there exists a quotient $\U$ as above. The diagram \eqref{eqn:diagram-UZ} implies, via the two properties we just discussed, that the intersection of $X(\Q_p)$ with the global Selmer variety is finite. Since this set contains $X(\Q)$, finiteness of the latter set follows. 

\par In fact, one can refine the above discussion by constructing a quotient which is an extension of $V$ by the direct sum of characters $\Q_p(1)^{\oplus (\mathrm{rk}_{\mathrm{NS}}-1)}$, resulting via the same reasoning in the following finiteness statement, which is a special case of Balakrishnan--Dogra \cite[Lemma 3.2]{BD18}.
\begin{theorem}\label{qclemma}
Suppose that $X$ is a smooth projective curve over $\Q$. Then $X(\Q)$ is finite whenever
\begin{equation}
r \ < \ g + \mathrm{rk}_{\mathrm{NS}} - 1.
\end{equation}
\end{theorem}
Many other instances of finiteness are known to follow from the method of Chabauty--Kim, and in general finiteness was proved by Kim \cite{Kim09} under the assumption of the Bloch--Kato conjecture. We will not discuss these results here, but rather turn to the question of how to explicitly determine the set $X(\Q)$. 

\section{The method of Chabauty--Kim: Effectivity}
\label{sec:CK2}
\par In this section, we discuss how to use the method of Chabauty--Kim to compute the
rational points on $X$ in the 
simplest instance of Theorem~\ref{qclemma}: the case $r=g$ and $\mathrm{rk}_{\mathrm{NS}} >1$. 
Whereas the method of Chabauty--Coleman relies on detecting global points via
\textit{linear}
relations in the image of $\per$, we will provide a
computable condition on filtered $\phi$-modules in the image of $\per$ to come
from a point in $X(\Q)$ via \textit{bilinear} relations, thereby addressing the third item in our wish
list in~\S\ref{subsec:quotients}.

\subsection{Heights on Selmer varieties}\label{sec:heights_selmer}
\par Looking for bilinear relations, one is naturally led to $p$-adic heights.
Classically, these were defined as bilinear pairings on $J(\Q)$
but since it is crucial that the non-abelian method of Chabauty--Kim factors through a non-abelian Selmer variety rather
than the abelian variety $J$, we instead prefer to utilize a more general approach due to
Nekov\'a\v r \cite[\S2]{Nek93}. Namely, he constructs a continuous bilinear pairing
\begin{equation}\label{eqn:NekovarGlobal}
h\colon \Hf(G_{\Q},V) \times \Hf(G_{\Q},V^*(1)) \lra \Q_p,
\end{equation}
depending on some auxiliary choices, including the choice 
of a splitting 
of the Hodge filtration 
\begin{equation}\label{eqn:splitting}
s\colon V_{\dR}/\Fil^0 V_{\dR} \lra V_{\dR}\,.
\end{equation}

\par 
The global height $h$ decomposes as a sum of local heights $h_v$, where $v$ runs through
the finite primes of $\Q$.
Briefly, the idea is to lift a pair in
$\Hf(G_{\Q},V) \times \Hf(G_{\Q},V^*(1))$ to a mixed extension of $p$-adic Galois
representations with graded pieces $\Q_p, V$ and $\Q_p(1)$ and to define $h_v$ on it.  As explained in~\cite[Section~5]{BD18}, we can construct such a representation from 
a torsor $P\in \Hf(G_{\Q}, \U)$, where $\U$ is attached to a N\'eron--Severi class as in \S\ref{subsec:correspondences}, by twisting a certain quotient of the
universal enveloping algebra of the $\Q_p$-unipotent \'etale fundamental group by $P$.
There is an analogous local construction for $P\in \Hf(G_{\Q_v}, \U)$.

\par 
We will assume throughout that $r=g$ and that the $p$-adic closure of $J(\Q)$ has finite
index in $J(\Q_p)$. \footnote{If the latter condition fails, we may apply classical Chabauty
as in~\S\ref{sec:rat_pts_expl}.} Then there are isomorphisms
\[
  \Hf(G_{\Q},V) \xrightarrow{\ \ \res_p \ \ }{} \Hf(G_{\Q_p},V) \xrightarrow{\ \ \log \ \
  }{} \HH^0(X_{\Q_p}, \Omega^1)^{\vee}.
\]
By Poincar\'e duality we obtain maps 
\[
  \pi\colon \Hf(G_{K}, \U) \lra \HH^0(X_{\Q_p},\Omega^1)^{\vee} \otimes \HH^0(X_{\Q_p},\Omega^1)^{\vee}
\]
for $K \in \{\Q,\Q_p\}$.

For ease of exposition, we shall assume for all $v\ne p$ that $h_v=0$ for torsors coming
from $X(\Q_v)$. 
The local height $h_p$ will be discussed in more detail below. 
The main point is that it
factors through $\D_{\cris}$, so we obtain the following
refinement of diagram~\eqref{eqn:diagram_kim}:
\begin{equation}\label{eqn:big_diagram}
\begin{tikzpicture}[->,>=stealth',baseline=(current  bounding  box.center)]
 \node[] (Hf) {$X(\Q )$};
  \node[below of=Hf,node distance=1.5cm] (Hfp) {$\Hf(G_{\Q},\U)$};
 \node[right of=Hf, node distance=3.4cm]  (Hff) {$X(\Q_p)$};
  \node[below of=Hff,node distance=1.5cm] (Hfpf) {$\Hf (G_{\Q_p},\U)$};
 \node[below of=Hfpf,node distance=1.5cm](Dieu) {$\HH^0(X_{\Q_p},\Omega^1)^{\vee} \otimes \HH^0(X_{\Q_p},\Omega^1)^{\vee}$};
 \node[right of=Hfpf,node distance=3.5cm](MF) {$\mathrm{MF}^{\phi}$};
 \node[below of=MF,node distance=1.5cm](Dieuf) {$\Q_p$\,,};
 \path (Hf) edge node[left]{$\scriptstyle \rho$ }(Hfp);
  \path (Hf) edge node[left]{}(Hff);
\path (Hfp) edge node[above]{ $\scriptstyle \res$ } (Hfpf);
\path (Dieu) edge node[below]{$\scriptstyle h$ } (Dieuf) ;
 \path (Hfp) edge node[below]{$\scriptstyle \pi$}(Dieu);
 \path (Hfpf) edge node[left]{$\scriptstyle \pi$}(Dieu);
 \path (Hff) edge node[left]{$\scriptstyle \rho$ } (Hfpf) ;
\path (MF) edge node[left]{$\scriptstyle h_p$} (Dieuf);
 \path (Hfpf) edge node[above]{\footnotesize $\D_{\cris}$} (MF);
 \path (Xp)  edge[dashed, bend left=20] node[above right]{\footnotesize $\per$} (MF);
\end{tikzpicture}
\end{equation}
where $h_p$ is now defined on the image of $\Hf (G_{\Q_p},\U)$.

If $(\psi_i)$ is a basis of the dual space of
$\HH^0(X_{\Q_p},\Omega^1)^{\vee} \otimes \HH^0(X_{\Q_p},\Omega^1)^{\vee}$, then there are
constants $\alpha_i\in \Q_p$ such that
$h = \sum_i \alpha_i\psi_i$. We deduce that the locally analytic function
  \begin{equation}\label{eq:qc_function}
Q\,\colon\, X(\Q_p) \lra \Q_p\,;\qquad x \mapsto \sum_i \alpha_i\psi_i(\pi(\rho(x))) -
 h_p(\per(x))  
  \end{equation}
vanishes along $X(\Q)$; furthermore, one can show
that it has only finitely many zeroes (see~\cite{BD19}). We can use this function for the explicit computation
of $X(\Q)$ if we have algorithms to
\begin{enumerate}
  \item[\upshape (i)] compute the $\alpha_i$ for a suitable explicitly computable basis $\psi_i$. 
  \item[\upshape (ii)] expand the function $x \mapsto h_p(\per(x))$ into convergent power series on
    residue disks.
\end{enumerate}
We can easily solve (i) given $x_1 ,\ldots ,x_m\in X(\Q)$ such that 
\[
  \{ \pi(\rho(x_i)) \}_{i = 1,\ldots, m}\ \ \mbox{is a basis for}\ \
  \HH^0(X_{\Q_p},\Omega^1)^{\vee} \otimes \HH^0(X_{\Q_p},\Omega^1)^{\vee}; 
\]
in this case we only need to compute $h_p(\per(x_i))$ and $\pi(\rho(x_i))$. 
If we choose an $\mathrm{End}_0(J)$-equivariant splitting in~\eqref{eqn:splitting}, then the
global height is also $\mathrm{End}_0(J)$-equivariant, thus reducing the number of points $x_i$ required.
Nevertheless, there need not exist enough points $x_i$, in which case we can solve (i)
using generators of $J(\Q)\otimes\Q$ and a construction of
$p$-adic heights on $J$ due to Coleman and Gross~\cite{CG89}.

\begin{remark}It is possible to write down functions vanishing in $X(\Q)$ with
finitely many zeroes when $r \ < \ g + \mathrm{rk}_{\mathrm{NS}} - 1$ using $p$-adic
heights~\cite[Proposition~5.9]{BD18}. More generally, one can extend Nekov\'a\v r's
construction to construct such functions when $r<g^2$, conditional on the conjecture of
Bloch--Kato, see~\cite[\S4]{BD19}. This has only been made explicit in the special
case of the Kulesz--Matera--Schost family of bielliptic genus 2
curves, see the (unconditional) Theorem~1.2 of \cite{BD19}.\end{remark}

\subsection{Local heights}
In the remainder of this section we focus on (ii). We first discuss in more detail the
local height $h_p$, following~\cite[\S4]{Nek93}. Let $P\in \Hf(G_{\Q_v}, \U)$ and denote
by $M_P$ the mixed extension of $\pi(P)$ mentioned above.
Then $h_p(M_P)$ is constructed using $\D_{\cris}(M_P)$, which is a mixed extension of filtered
$\phi$-modules with graded pieces $\Q_p, V_{\dR} \colonequals
\HH^1_{\dR}(X_{\Q_p})^{\vee} = \D_{\cris}(V)$ and $\Q_p(1)$.

For simplicity, we only describe $h_p$ on the image of $X(\Q_p)$.
The family $(\D_{\cris}(M_{\rho(x)}))_{x}$ interpolates in the following sense:
There is a filtered connection $\cA_Z=\cA_Z(b)$ with Frobenius structure such that we have 
\begin{equation}\label{}
  \D_{\cris}(M_{\rho(x)}) \simeq  x^*\cA_Z\quad\text{for all}\;\; x\in X(\Q_p).
\end{equation}

Suppose that we have isomorphisms \[
\begin{array}{lllll}
  s^\phi(b,x) &\colon& \ \Q_p \oplus V_{\dR} \oplus \Q_p(1) & \stackrel{\sim }{\lra } & x^*\cA_Z \\
  s^{\Fil}(b,x) &\colon& \ \Q_p \oplus V_{\dR} \oplus \Q_p(1) & \stackrel{\sim }{\lra } & x^*\cA_Z \\
\end{array}
\]
where $s^{\phi}$ is Frobenius-equivariant, and $s^{\Fil}$ respects the filtrations, and
suppose we can write them as 
\begin{equation}\label{eqn:definitions-abcs}
s^\phi(b,x) =\left(
\begin{array}{ccc} 
1 & 0 & 0 \\
\boldsymbol{\alpha }_{\phi}(b,x) & 1 & 0 \\
\gamma_{\phi}(b,x) & \boldsymbol{\beta }^{\intercal}_{\phi}(b,x) & 1 \\
\end{array}
\right) \qquad 
s^{\Fil}(b,x) = \left( 
\begin{array}{ccc}
1 & 0 & 0 \\
0 & 1 & 0 \\
\gamma_{\Fil}(b,x) & \boldsymbol{\beta }^{\intercal}_{\Fil} (b) & 1 \\
\end{array}
\right).
\end{equation} Note that we make a choice of basis differentials on the affine open $Y$
(see~\S~\ref{sec:hodge}) so that $s^\phi(b,x)$ and $s^{\Fil}(b,x)$ are of this form.
The splitting $s$ in~\eqref{eqn:splitting} induces idempotents 
\[
\begin{array}{ccccc}
  s_1& \colon& V_{\dR} & \lra & s(V_{\dR}/\Fil^0 V_{\dR}) \\
  s_2& \colon& V_{\dR} & \lra & \Fil^0 V_{\dR}. 
\end{array}
\]
With respect to our choices, Nekov\'a\v r's local height at $p$ is
\begin{equation}\label{eq:ht_formula}
  h_p(\per(x))=
 \gamma_{\phi}(b,x) 
 - \gamma_{\Fil}(b,x) 
 -\boldsymbol{\beta}^\intercal_{\phi}(b,x) \cdot s_1 (\boldsymbol{\alpha}_{\phi}) (b,x) 
 - \boldsymbol{\beta }^\intercal_{\Fil}(b) \cdot s_2 (\boldsymbol{\alpha}_{\phi})(b,x).
  \end{equation}

So in order to solve (ii) we need to compute the entries of~\eqref{eqn:definitions-abcs}, which means
computing the Hodge filtration and the Frobenius structure on $\cA_Z$. 
For (i), we also need to explicitly compute the composition $\pi \circ \rho$. 
With respect to the dual basis of our chosen basis differentials on $Y$, the map $\pi \circ \rho$ is given by 
\begin{align}\label{eq:mixed_cpts}
\pi \circ \rho: Y(\Q_{p}) & \to H^0(X_{\Q_p}, \Omega^1)^{\vee} \otimes  H^0(X_{\Q_p}, \Omega^1)^{\vee}\,\\
x & \mapsto \boldsymbol{\alpha}_{\phi}(b,x)^{\intercal}\cdot \left(\begin{array}{c} I_g
 \nonumber \\
  \hline 0_g \end{array} \right) \otimes  (\boldsymbol{\beta}_{\phi}^{\intercal}(b,x) -
  \boldsymbol{\beta}_{\Fil}^{\intercal}(b))\cdot \left(\begin{array}{c} 0_g \\ \hline I_g \end{array} \right). 
\end{align}
Note in particular that the first factor is the Abel-Jacobi map $\AJb(x)$, sending $x$ to
the functional $\omega \mapsto \int_{b}^{x} \omega$. 

\subsection{Computing the Hodge filtration}\label{sec:hodge}

We work in an affine open subset $Y$ of $X$.  Suppose that we have $\#(X\setminus Y)(\overline{\Q}) = d$ and choose differentials $\omega_0, \ldots, \omega_{2g+d-2} \in H^0(Y_{\overline{\Q}},\Omega^1)$ on $Y$ such that the following conditions are satisfied:
\begin{enumerate}\item The differentials $\omega_0, \ldots \omega_{2g-1}$ are of the second kind (residue zero) on $X$  and form a symplectic basis of $H^1_{\dR}(X_{\Q})$ with respect to the cup product pairing. We let $\boldsymbol{\omega}$ denote the column vector $(\omega_0, \ldots \omega_{2g-1})^{\intercal}$.
\item The differentials $\omega_{2g}, \ldots, \omega_{2g+d-2}$ are of the third kind (all poles have order one) on $X$.
\end{enumerate}

Universal properties give that the rank $2g+2$ vector bundle $\mathcal{A}_Z$ has a
connection, a Hodge filtration, and a Frobenius structure, as discussed in \cite[\S 4,5]{BDMTV19}. Here, we give algorithms that describe these objects.

Recall that we have a non-trivial class $Z$ in $\mathrm{NS}(J)$ mapping to~0 in
$\mathrm{NS}(X)$.  This is equivalent to the choice of an endomorphism of $\HH^1_{\dR}(X)$ satisfying several conditions (see \cite[\S 4.4]{BDMTV19}), and we describe a method to compute this in the case of modular curves in Section \ref{sec:Examples}.  We denote the matrix of the correspondence $Z$ on $\HH^1_{\dR}(X/\Q)$ also by $Z$, where we act on column vectors.

Choose a trivialization $$s_0\, \colon\, \mathcal{O}_Y \otimes (\Q_p \oplus V_{\dR}
\oplus \Q_p(1))\, \rightarrow\, \mathcal{A}_Z|_Y$$ such that, with respect to this trivialization, the connection $\nabla$ on $\mathcal{A}_X$ is given by
$$\nabla = d + \Lambda,$$ where $$\Lambda =  -\left( \begin{array}{ccc}
0 & 0 & 0 \\
\boldsymbol{\omega} & 0 & 0 \\
\eta & \boldsymbol{\omega}^{\intercal}Z & 0 \\
\end{array} \right),$$ where $\eta$ is a differential of the third kind on $X$ that is uniquely determined by the following two properties:
\begin{enumerate}\item It is in the space spanned by $\omega_{2g}, \ldots, \omega_{2g+d-2}$, and 
\item The connection $\nabla$ extends to a holomorphic connection on all of $X$.
\end{enumerate}
The Hodge filtration on $\mathcal{A}_Z$ is determined completely from the Hodge filtration on its graded pieces, via universal properties. Here is an algorithm to compute the Hodge filtration:

\begin{algorithm}[Computing the Hodge filtration on $\mathcal{A}_Z$]\label{algo:hodge} \;$\quad$
\begin{enumerate}
\item Let $L/\Q$ denote a finite extension over which all the points of $X \setminus Y$ are defined. Compute local coordinates at each $x \in (X \backslash Y)(L)$.
\item For each $x \in (X \backslash Y)(L)$, compute power series for $\boldsymbol{\omega}_x$, the expansion of the vector of differentials $\boldsymbol{\omega}$ at $x$ to large enough precision, which means at least mod ${t_x^{d_x}}$, where $d_x$ is the order of the largest pole occurring.
\item Compute the vector $\boldsymbol{\Omega}_x$, defined by $$d\boldsymbol{\Omega}_x = - \boldsymbol{\omega}_x.$$
\item Compute $\eta$ as the unique linear combination of $\omega_{2g}, \ldots , \omega_{2g+d-2}$ such that 
\[
d\boldsymbol{\Omega}_x^{\intercal}Z\boldsymbol{\Omega}_x -  \eta
\]
has  residue zero at all $x \in (X \backslash Y)(L)$. To do this, carry out the following:
\begin{enumerate} \item Using local coordinates at each $x \in (X \backslash Y)(L)$, rewrite $\omega_{2g}, \ldots , \omega_{2g+d-2}$. 
\item Solve for $\eta$ by comparing residues.
\end{enumerate}
\item Solve the system of equations for $g_x$ in $L (\! (t_x )\! )/L [\! [t_x ]\! ]$ such that $$dg_x =  \boldsymbol{\Omega}_x^{\intercal}Z d\boldsymbol{\Omega}_x -  \eta.$$ 
\item Compute the vector of constants $\mathbf{b}_{\Fil } = (b_g,\ldots, b_{2g-1})\in \Q^g$ and the function $\gamma_{\Fil}$
  characterized by $\gamma_{\Fil}(b)=0$ and 
\begin{align} \label{eqn:final-hodge-equation}
g_x +\gamma _{\Fil } - \mathbf{b}_{\Fil }^{\intercal}N^{\intercal}\boldsymbol{\Omega }_x   -\boldsymbol{\Omega}_x ^{\intercal }Z NN^{\intercal}\boldsymbol{\Omega}_x \in L[\! [t_x ]\! ]
\end{align}
where $N$ is the $2g\times g$ matrix which has the zero matrix  of dimension $g$ and the
    identity matrix of dimension $g$ as blocks. Set $\boldsymbol{\beta}_{\Fil}=
    \boldsymbol{\beta}_{\Fil}(b) = (0,\ldots, 0, b_g, \ldots, b_{2g-1})^{\intercal}$.
\end{enumerate}
\end{algorithm}

\begin{remark}\label{hypeta}We note that \cite[Lemma 6.5]{BD19} simplifies some of the calculations in the case of a hyperelliptic curve $X$: in this case, we have that $\eta=0$ and $\boldsymbol{\beta}_{\Fil} = (0,\ldots,0)^{\intercal}$.\end{remark}

\subsection{Computing the Frobenius structure}
The Frobenius structure on $\cA_Z$ can be determined explicitly in terms of double Coleman integrals, as discussed in \cite[\S 5]{BDMTV19}. Here is an algorithm to compute it:

\begin{algorithm}[Computing the Frobenius structure on $\mathcal{A}_Z$]\label{algo:frob} \;$\quad$
\begin{enumerate}
\item Use Tuitman's algorithm \cite{Tui16, Tui17} to compute the matrix of Frobenius $F$
  and a vector $\boldsymbol{f}$ of overconvergent functions such that $$\phi^*\boldsymbol{\omega} = d\boldsymbol{f} + F\boldsymbol{\omega},$$
    where $\phi$ is a certain lift of Frobenius.
\item Let $b_0,x_0$ be Teichm\"uller representatives of $b,x$ respectively. Compute the matrix 
\[
\mathrm{A} = \mathrm{I}(x,x_0)^+ \cdot \mathrm{I}(b_0,b)^-,
\] 
where we define for any pair $x_1,x_2 \in X(\Q_p)$ the parallel transport matrices
\[
\mathrm{I}^{\pm}(x_1,x_2) = 
\left( \begin{array}{ccc}
1 & 0 & 0 \\
\int_{x_1}^{x_2} \boldsymbol{\omega} & 1 & 0 \\
\int_{x_1}^{x_2} \eta + \int_{x_1}^{x_2} \boldsymbol{\omega}^{\intercal}Z\boldsymbol{\omega} & \pm \int_{x_1}^{x_2} \boldsymbol{\omega}^{\intercal}Z & 1 \\
\end{array} \right), 
\]
where $\eta$ is as computed in Algorithm \ref{algo:hodge} (see also Remark \ref{hypeta}).

\item Explicitly solve the system 
\[
\left\{ 
\begin{array}{lll}
d\mathbf{g}^\intercal & = & d\mathbf{f}^\intercal ZF, \\
dh & = & \boldsymbol{\omega}^\intercal F^\intercal Z\mathbf{f} + d\mathbf{f}^\intercal Z\mathbf{f}-\mathbf{g}^\intercal \boldsymbol{\omega}
+\phi^*\eta - p\eta, \\
h(b) & = & 0. \\
\end{array} 
\right.
\]
Then compute the matrix
\begin{equation*}\label{eqn:frob-triv-Teich}
\mathrm{M}(b_0,x_0) = \left( \begin{array}{ccc}
1 & 0 & 0 \\
(I-F)^{-1}\mathbf{f} & 1 & 0\\
\frac{1}{1-p}\left(\mathbf{g}^{\intercal}(I-F)^{-1}\mathbf{f} +h\right) & \mathbf{g}^{\intercal}(F-p)^{-1} & 1 \\
\end{array} \right)(x_0).
\end{equation*}

\item Finally, compute the matrix 
\begin{equation*}\label{eqn:frob-triv-non-Teich}
 s_0 ^{-1}(b,x) \circ s^\phi (b,x) = \mathrm{A} \cdot \mathrm{M}(b_0,x_0) = \left( \begin{array}{ccc}
1 & 0 & 0 \\
\boldsymbol{\alpha}_{\phi}(b,x) & 1 & 0\\
\gamma_{\phi} (b,x) & \boldsymbol{\beta}^\intercal_{\phi}(b,x) & 1 \\
\end{array} \right).
\end{equation*}
\end{enumerate}
\end{algorithm}

\begin{remark}\label{hypfrob} If $X$ is a hyperelliptic curve, say the smooth projective
  model of the affine curve $Y\colon y^2=f(x)$, where $f$ is monic and has no repeated roots, then we can use
  Kedlaya's algorithm~\cite{Ked01} or Harrison's generalization \cite{harrison} in Step~(1) above. In fact, Tuitman's approach generalizes the approach of
  Kedlaya and Harrison.  Note that the \texttt{SageMath} implementation of Kedlaya's algorithm takes the convention that Frobenius acts on columns, while the \texttt{Magma} implementation of Tuitman's algorithm as used here takes the convention that Frobenius acts on rows and thus differs by a transpose.
\end{remark}
\begin{remark}Computing the action of Frobenius in Step (1) gives us a way to compute Coleman integrals: in particular, if $b_0 = \phi(b_0)$ and $x_0 = \phi(x_0)$ are Teichm\"uller points, we compute the Coleman integral as
\[
\int_{b_0}^{x_0} \boldsymbol{\omega} = (1 - F)^{-1} \left(\boldsymbol{f}(x_0) - \boldsymbol{f}(b_0)\right)\,.
\]
\end{remark}

\section{Examples}
\label{sec:Examples}

We illustrate the practicality of the method of Chabauty--Kim discussed in Section~\ref{sec:CK2} by applying it to three new examples of curves whose rational points were previously unknown. They are all curves of the form 
\begin{equation}
X_0(N)^+ \colonequals X_0(N)/w_N
\end{equation}
where $N$ is prime and $w_N$ is the Atkin--Lehner involution, and therefore they have a
unique rational cusp, and their non-cuspidal rational points classify unordered pairs of elliptic curves that are related by an $N$-isogeny. We consider the cases $N=67$, $73$, and $103$. For each value of $N$, the curve $X_0(N)^+$ is of genus $2$ and its Jacobian has real multiplication. Thus, the rank of the N\'eron--Severi group is equal to $2$, and the method outlined in Section~\ref{sec:CK2} produces exactly one non-trivial locally analytic function on $X_0(N)^+(\Q_p)$ vanishing on the set of rational points $X_0(N)^+(\Q)$. Hence, unlike in the Chabauty--Coleman example at the end of Section \ref{sec:Intro}, we need in addition the Mordell--Weil sieve (see \S\ref{subsec:MordellWeil}) to extract the set of rational points from the larger quadratic Chabauty set. 

\par We discuss the computation for $N=67$ in some detail and briefly summarize the cases
$N=73$ and $N = 103$. These computations use the computer algebra system {\tt
Magma}~\cite{BCP97} and were done by Best, Bianchi, and Triantafillou, mostly at
the workshop ``Arithmetic Statistics and Diophantine Stability'' at the Fondation des
Treilles in July 2018. 

\begin{remark} In~\cite{BGX19}, the authors apply a combination of elliptic curve Chabauty with
covering techniques to determine the rational points on $X_0(N)^+$ for several
composite squarefree values of $N$ such that $X_0(N)+$ has genus~2. It would be interesting to
determine the rational points on the 13 remaining hyperelliptic curves $X_0(N)+$ for
squarefree $N$; all of them have genus~2.
\end{remark}

\subsection{An explicit model for $X_0(67)^+$. }
As is explained in \cite{Murabayashi,Gal96}, an affine model for the genus $2$ curve $X_0(67)^+$ can be found explicitly as follows. Let $f$ be the unique, up to conjugation, newform of level $67$, weight $2$, which is furthermore invariant under the Atkin--Lehner involution $w_{67}$. The complex vector space spanned by $f$ and its Galois conjugate $f^c$ is isomorphic to the space of regular differentials on $X_0(67)^+$, and we may choose a basis $g_1$ and $g_2$ for this space such that $g_1=q+\cdots$ and $g_2=q^2+\cdots$. Note that $f$ and $f^c$ can be computed up to arbitrary $q$-adic precision using {\tt Magma} \cite{BCP97}. Then $x=\frac{g_1}{g_2}$ and $y=\frac{q}{g_2}\frac{dx}{dq}$ are related by an equation of the form $y^2=p(x)$, for some monic polynomial $p(x)$ of degree $6$ whose coefficients can be determined from the $q$-expansions. Such an equation gives a model for $X_0(67)^+$; while $g_2$ is unique, a certain choice of $g_1$ yields
\[
Y\colon y^2 =  x^6 + 2 x^5 + x^4 - 2x^3 + 2x^2 - 4x + 1 \hspace{2em} [\  \equalscolon f_{67}(x)\ ]\,.
\]
See \cite{Murabayashi} for more details; for other examples of computations of models of higher genus modular curves, see \cite{Gal96}.
The projective closure $X$ adds two points at infinity, $\infty^+$ and $\infty^-$, corresponding to $(1:1:0)$ and $(1:-1:0)$ respectively. By an explicit search, we quickly find several points in $X(\Q)$. Indeed, 
\begin{equation}\label{eqn:knownpts67}
X(\Q) \supset \{\infty^+, \infty^-, (0,\pm 1), (-1, \pm 3), (1, \pm 1), (-2, \pm 7)\}\,.
\end{equation}

Our goal is to use the machinery set up in Section \ref{sec:CK2}, combined with the Mordell--Weil sieve, to show that $X(\Q)$ consists precisely of these $10$ points.

\par Using the explicit model $Y$, several arithmetic properties of $X_0(67)^+$ can be deduced. For instance, {\tt Magma}'s implementation of 2-descent shows that
that the rank\footnote{Since the Jacobian $J_0(N)^+$ has real multiplication over $\Q$ for every prime $N>2$, its Mordell-Weil rank
over $\Q$ is necessarily a multiple of the genus.} of $J_0(N)^+(\Q)$ is exactly $2$. Alternatively, one can avoid the use of a model and draw the same conclusion from the Gross--Zagier--Kolyvagin--Logachev theorem \cite{GZ86, KL89}, by computing that (provably \cite[Chapter 3]{Ste00}) $L(f,1) = 0$ and (numerically \cite{Cre97, Dok04}) $L'(f,1)\neq 0$.

\subsection{The reduction of $X_0(N)^+$. }\label{sec:semistable_model} 
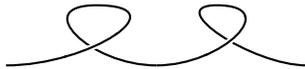
\begin{figure}
  \centering
    \begin{tikzpicture}[scale = 0.4]
    \begin{knot}
      \strand[thick] (-5,0) to [out=right, in=right, looseness=2] (-2,2);
      \strand[thick] (-2,2) to [out=left, in=left, looseness=2] (0,0);
      \strand[thick] (0, 0) to [out=right, in=right, looseness=2] (2,2);
     \strand[thick] (2,2) to [out=left, in=left, looseness=1.3] (5,0);
    \end{knot}
  \end{tikzpicture}
      \caption{The reduction of $\mathcal{X}_0(N)^+$ at $N$. }\label{semist}
  \end{figure}

Recall that the
method outlined in Section \ref{sec:CK2} uses some global and local $p$-adic heights in
the sense of Nekov\'a\v r. Although these depend on some auxiliary choices that we have
not made yet at this stage, we have already remarked that we can always ignore all the
local heights at primes $v\neq p$ of potential good reduction. More generally, by work of
Betts--Dogra \cite[Corollary 1.2.2]{BeD19}, the map $X(\Q_v)\to \Q_p$ induced by the local
height at $v\neq p$ takes at most as many values as the number of irreducible components
of a regular semi-stable model at $v$.  
Note that $X_0(N)^+$ has good reduction at all primes away from $N$.  Using an
argument analogous to~\cite[Theorem~6.6]{BDMTV19}, we can 
show that for all primes $N >2$ there is a regular semi-stable model $\mathcal{X}_0(N)^+$
of $X_0(N)^+$ whose special fibre is isomorphic to a projective line intersecting itself $g$ times, where $g$ is the genus of $X_0(N)^+$ (see Figure~\ref{semist}). The self-intersections correspond to conjugate pairs of supersingular $j$-invariants in $\F_{N^2}\setminus \F_{N}$ (see \cite[V, \S1]{DR73} and \cite[\S3]{Ogg75}).
In particular, the special fibre of $\mathcal{X}_0(N)^+$ consists of only one component, so the work of Betts--Dogra implies that there are no non-trivial contributions at $v\neq p$. 

\subsection{Preliminary choices}
 \subsubsection*{A prime $p$ and a base point $b$} Since by \S\ref{sec:semistable_model} the curve $X_0(N)^+$ has good reduction at all primes away from $N$, we could let our fixed $p$ be any prime different from $N$; we pick $p=11$. This choice may seem slightly peculiar to the reader familiar with the classical Chabauty--Coleman method, where it is often advantageous to choose the smallest possible prime of good reduction. The prime $11$ has two main advantages for our purposes. First, the polynomial $f_{67}$ has no linear factors over $\Q_{11}$ and, as a result, the lift of Frobenius that we use in \S\ref{sec:frob_structure}
extends to all of $X(\Q_{11})$. 
 While it is possible to deal with disks containing a point with zero $y$-coordinate by  working with a different lift of Frobenius or by using the trick discussed in \cite[\S 5.5]{BDMTV19}, our choice of $p$ makes both the exposition and the computation significantly shorter. The second advantage of the prime $11$ is somewhat post-hoc, coming from the final Mordell-Weil sieve step. Indeed, it turns out that $J_0(N)^+(\F_q)$ has order divisible by $11^2$ for several small primes $q$ (including $q = 31$ and $q = 137$), which makes the Mordell-Weil sieve particularly efficient for proving that points of $X(\Q_{11})$ are \emph{not} in $X(\Q)$.

 	 We choose $b = (1,1)$ for the base point. Note that $b$ lies in both our affine patch and in the affine patch at infinity. One advantage of $b$ over other possible base points is that $b$ will be a Teichm\"uller point for a convenient lift of Frobenius.

\subsubsection*{A basis for the de Rham cohomology of $X_0(67)^+$}
\label{subsec:symplectic-basis}
It is well known that $\HH^0(Y_{\overline{\Q}},\Omega^1)$ has basis given by (the classes of) the differentials
\begin{equation}
\label{eq:standard_basis}
\left\{\frac{\, dx}{y}, \frac{x\, dx}{y}, \frac{x^2 \, dx}{y}, \frac{x^3 \, dx}{y}, \frac{x^4 \, dx}{y}
\right\}\,
\end{equation}
and that, inside of $\HH^0(Y_{\overline{\Q}},\Omega^1)$, we can identify $\HH^1_{\dR}(X)$ with those differentials which have residue $0$ at both points of $X\smallsetminus Y = \{\infty^+, \infty^-\}$. By working with the expansion of each differential in \eqref{eq:standard_basis} in terms of the uniformizer $t_{\infty^{\pm}}=x^{-1}$ at $\infty^{\pm}$, we construct a new basis $\omega_0,\dots, \omega_4$ satisfying the properties (1) and (2) of \S\ref{sec:hodge}. In particular, we may take
\begin{align*}
\omega_0 = -\frac{\, dx}{y},\quad \omega_1 =  (-1 - x) \cdot \frac{\, dx}{y},\quad \omega_2 =  (-2 + x - x^3 - x^4) \cdot \frac{\, dx}{y},\\
 \omega_3=\frac{1}{2}\left(1 -  x^2 -  x^3\right) \cdot \frac{\, dx}{y},\quad \omega_4 = (-x - x^2) \cdot \frac{\, dx}{y}\,.
\end{align*}

From now on, $\boldsymbol{\omega}$ will denote the column vector $(\omega_0,\dots,\omega_3)^{\intercal}$.

\subsubsection{A N\'eron--Severi class}\label{sec:neronseveri}
The choice of a N\'eron--Severi class $Z$ as in Section \ref{sec:CK1} is equivalent to the choice of an endomorphism of $\HH^1_{\dR}(X)$, satisfying a list of conditions (see \cite[\S4.4]{BDMTV19}). Let $\ell$ be a prime of good reduction for $X$. In order to compute the action of the Hecke operator $T_{\ell}\in \mathrm{End}(\HH^1_{\dR}(X/\Q_{\ell}))$ on the whole of $\HH_{\dR}^1(X/\Q_{\ell})$, rather than just on $\Fil^0 \HH^1_{\dR}(X/\Q_{\ell})$, we use the Eichler--Shimura formula
\[
T_{\ell} = \Fr_{\ell}^{\intercal} + \ell \cdot (\Fr_{\ell}^{\intercal})^{-1}\,.
\]
The matrix of Frobenius $\Fr_{\ell}$ with respect to the basis $\boldsymbol{\omega}$ may be computed using Tuitman's algorithm (we briefly postpone a discussion of this to Step (1) of \S\ref{sec:frob_structure}, since this matrix for $\ell = p$, as well as one additional output of Tuitman's algorithm, are both needed at that step), and we identify the operator $T_{\ell}$ with its matrix representation with respect to $\boldsymbol{\omega}$. Note that the Eichler--Shimura formula holds for $X_0(67)$ and thus for $X_0(67)^+$, since the Atkin--Lehner involution commutes with $T_{\ell}$ at all $\ell\neq 67$.

To obtain from $T_{\ell}$ an endomorphism corresponding to a class $Z\in \mathrm{NS}(J)$ which maps to zero in $\mathrm{NS}(X)$, we first consider $\Tr(T_{\ell})\cdot I_{4}- 4T_{\ell}$, which has trace zero, and then multiply on the right by the inverse of the cup product matrix on $\boldsymbol{\omega}$. For example, choosing $\ell = 11$, we obtain the non-trivial endomorphism with matrix representation
\[
Z = \left( \begin{array}{cccc}
   0 & -8 & 12 &   8 \\
   8 &  0 & -8 & -12 \\
 -12 &  8 &  0 &   0 \\
  -8 & 12 &  0 &   0 
\end{array} \right)\,.
\]

Since the N\'eron--Severi group has rank $2$, choosing a different Hecke operator would only change the matrix $Z$ by a multiplicative constant factor.

\begin{remark}
Using Tuitman's algorithm, we can compute the entries of $T_{\ell}$ only up to some $\ell$-adic precision. In our case, this would suffice to carry out the steps of the quadratic Chabauty computation, since we have chosen $\ell = p$. It should however be possible to prove that $Z$ is given exactly by the above matrix, and we may assume that this is the case, as doing so does not affect the computation in any crucial way.
\end{remark}

\subsection{Hodge Filtration on $\cA_{Z}$}
\label{sec:hodge_fil}
We now compute the Hodge filtration of the vector bundle $\cA_{Z}$ attached to our choice of N\'eron--Severi class $Z$ and base point $b$. Since the curve $X$ is hyperelliptic, by Remark \ref{hypeta} we only need to compute $\gamma_{\Fil}$, and we can do so using a simplified version of Algorithm \ref{algo:hodge}. In particular, for each point at infinity $\infty^{\pm}$, we can compute $\boldsymbol{\Omega}_{\infty^{\pm}}$ and $g_{\infty^{\pm}}$ by formal integration of Laurent series in the uniformizer $t_{\infty^{\pm}}$. Following the steps, we then find that $\gamma_{\Fil}$ has a pole of exact order $1$ at $\infty^{\pm}$ with residue $-8$. Since $\gamma_{\Fil}$ must vanish at $b$, we conclude that
\begin{equation*}
\gamma_{\Fil} = -8x+8.
\end{equation*}

\subsection{Frobenius structure on $\cA_{Z}$}
\label{sec:frob_structure}
We compute the Frobenius structure on $\cA_{Z}$ using Algorithm \ref{algo:frob}.

{\bf Step (1):} We first fix a lift of Frobenius $\phi$. We take $\phi(x) = x^p$, and extend to $\Q_{11}[x]$ by linearity. Since $f_{67}$ has no
zeros over $\F_{11}$, we extend this lift to a strict open neighborhood of the tube $] Y_{\F_{p}}[$\footnote{the tube consists of all points	reducing to $Y_{\F_{p}}$} 
by expanding
\[
\phi(y) = \sqrt{ \phi(f_{67}(x))} = y^p \cdot \left( 1 + \frac{\phi(f_{67}(x)) - f_{67}(x)^p}{y^{2p}} \right)^{1/2}\,.
\]
as an \emph{overconvergent} Laurent series in $\Q_p[\! [x, y, y^{-1}]\! ]$. This lift naturally extends to one-forms.

Next, we compute $p$-adic approximations of $F$ and $\boldsymbol{f}$ using Tuitman's algorithm  \cite{Tui16,Tui17}, a generalization of Kedlaya's algorithm which incorporates Lauder's fibration method \cite{Lau06}.
Roughly speaking, we first compute $\phi^* \omega_i$. Then, we reduce pole orders by iteratively subtracting differentials of overconvergent functions (constructed by solving linear systems) until $\phi^{*} \omega_i$ has been reduced to a cohomologous linear combination of basis differentials $\sum_{j} F_{ji} \omega_j$. The sum $f_i$ of the functions from each step satisfies
\[
\phi^{*} \omega_i = \sum_{j}F_{ij} \omega_j + df_i\,.
\]
Note that in our working example, this $F$ is the matrix $\mathrm{Fr}_{\ell}$ that was computed in \S\ref{sec:neronseveri} since $\ell$ was chosen there to be $p=11$ as well.
 
{\bf Step (2):} Since $b = (1,1)$ is a Teichm\"uller point for $\phi$, $I(b_0,b)^{-} = I(b,b)^{-}$ is an identity matrix. To compute the $I(x, x_0)^+$ on each residue disk, we expand the $\omega_{i}$ in terms of a uniformizer near each Teichm\"uller point $x_0$ and integrate formally. To compute $\int_{x_0}^{x}\boldsymbol{\omega}^{\intercal}Z\boldsymbol{\omega}$, we expand, formally integrate, multiply terms, and formally integrate again, as in steps (3) and (5) of Algorithm \ref{algo:hodge}.

{\bf Step (3):} The matrices $Z$ and $F$ are constants, so  $\boldsymbol{g}^{\intercal} = \boldsymbol{f}^{\intercal} Z F$.
We approximate $h$ by iteratively ``reducing'' a $p$-adic approximation $(dh)_{\sim}$ to  
$
dh = \boldsymbol{\omega}^{\intercal} F^{\intercal} Z \boldsymbol{f} + d\boldsymbol{f} Z \boldsymbol{f} - \boldsymbol{g}^{\intercal} \boldsymbol{\omega} + \phi^{*} \eta - p \eta
$
as in Tuitman's algorithm until we find $a_j \in \Q_{11}$ and an overconvergent function $h_{\sim}(x)$
so that
\[
(dh)_{\sim} = \sum_{j} a_j \omega_j + d(h_{\sim})\,.
\]
Then $h_{\sim}(x)-h_{\sim}(b)$ approximates $h(x)$.
The remainder of Steps (3) and (4) are straightforward. The terms $\boldsymbol{\alpha}_{\phi}(b,x), \boldsymbol{\beta}_{\phi}(b,x), \gamma_{\phi}(b,x)$ cannot be expressed compactly, so we omit them here.

\subsection{The local $p$-adic height and a finite set of $p$-adic points containing $X(\Q)$.}

We have now assembled all ingredients to compute the quadratic Chabauty function  from~\eqref{eq:qc_function}, whose finite set of zeroes contains $X(\Q)$. 
To find the constants $a_i$ in~\eqref{eq:qc_function}, we use the discussion at the end of \S\ref{sec:heights_selmer}. 

Set $K\colonequals \Q(\sqrt{5}) = \mathrm{End}_0(J_0(67)^+)$ and $K_{p} = K \otimes_{\Q} \Q_p$. If we pick a $K$-equivariant splitting $s$ of the Hodge filtration in formula~\eqref{eq:ht_formula}, then 
the global height $h$ factors through the tensor product $H^0(X_{\Q_p}, \Omega^1)^{\vee} \otimes_{K_{p}} H^0(X_{\Q_p}, \Omega^1)^{\vee}$.
We now choose auxiliary points $x_1 = (-2,7), x_2 = (-1,3) \in X(\Q)$. Since
$\AJb(x_1) = [\omega \mapsto \int_{b}^{x_1} \omega]$ is nonzero, $\AJb(x_1)$ is a $K_{p}$-basis for $H^0(X_{\Q_p},\Omega^1)^{\vee}$.
Using~\eqref{eq:mixed_cpts}, we compute
$(\pi \circ \rho)(x_i)$ in this basis. 

We compute $h(\pi(\rho(x_i)))   = h_p(\per(x_i))$ using~\eqref{eq:ht_formula}, the results
of \S\ref{sec:hodge_fil}, \S\ref{sec:frob_structure}
and the splitting $s$ associated to the $K$-equivariant basis
$(\omega_0,\omega_1,\omega_2,\omega_3-\omega_1)$.
Writing $\psi_1$ for the projection onto the ``rational part'' and $\psi_2$ for the projection onto the ``$\sqrt{5}$ part,'' we find that the function sending $x\in X(\Q_p)$ to
\begin{align}\label{eq:final_qc_fn}
Q(x) \colonequals h_p(\per(x)) 
& - ( 5 \cdot 11 + 2\cdot 11^2 + 5\cdot 11^3 +  0\cdot 11^4 + \cdots ) \cdot \psi_1(\pi(\rho(x)))\nonumber\\ 
& \qquad + (4 \cdot 11 + 0\cdot 11^2 + 4\cdot 11^3 + 0\cdot 11^4 + \cdots ) \cdot \psi_2(\pi(\rho(x))) 
\end{align}
vanishes for all $x\in X(\Q)$.

We expand $Q$ as a power series on each residue disk, find the roots, and repeat the
computation on an affine patch containing the points at infinity to find a finite subset
of $X(\Q_{11})$ which contains $X(\Q)$. Using a Newton polygon argument, we find that
every root of $Q$ is simple. In addition to the $10$ known rational points, we find $14$
additional $11$-adic zeros of $Q$ (listed in Table \ref{table:Qp-pointsX067+}). To show
that these points are not rational, we turn to the Mordell-Weil sieve, described in the
following subsection.

\begin{table}
\begin{tabular}{|c||r|} \hline
	Disks & $x$-coordinates of candidate points \\ \hline \hline
	$](0,\pm 1)[$  & $0$ \\
	           & $0 +  7 \cdot 11 +  0 \cdot 11^2 + 3 \cdot 11^3 + 3 \cdot 11^4 + \cdots$ \\ \hline 
	$](1,\pm 1)[$  & $1$ \\
               & $1 + 6 \cdot 11 + 6 \cdot 11^2 + 8 \cdot 11^3 + 7 \cdot 11^4 + \cdots$ \\  \hline 
	$](6,\pm 5)[$  & $6 + 5 \cdot 11 + 8 \cdot 11^2 + 2 \cdot 11^3 + 4 \cdot 11^4 + \cdots$ \\ 
               & $6 + 7 \cdot 11 + 0 \cdot 11^2 + 5 \cdot 11^3 + 1 \cdot 11^4 + \cdots$ \\ \hline 
	$](-2, \pm 7)[$ & $-2$ \\
               & $9 + 10 \cdot 11 + 1 \cdot 11^2 + 8 \cdot 11^3 + 0 \cdot 11^4 + \cdots$ \\ \hline	         
	$](-1,\pm 3)[$ & $-1$ \\ 
               & $10 + 3 \cdot 11 + 9 \cdot 11^2 + 10 \cdot 11^3 + 1 \cdot 11^4 + \cdots $\\ \hline
	$]\infty^{\pm}[$ & $\infty$ \\ 
& $2\cdot 11^{-1} + 4 + 10 \cdot 11 + 9 \cdot 11^2 + 8 \cdot 11^3 + 7 \cdot 11^4 + \cdots$ \\ \hline 
\end{tabular}
\caption{A set of $24$ points of $X_{0}(67)^+(\Q_{11})$ containing $X_{0}(67)^+(\Q)$.}
\label{table:Qp-pointsX067+}
\end{table}

\subsection{The Mordell--Weil sieve}
\label{subsec:MordellWeil}

We assume we are given a smooth projective curve $X/\Q$, $p$ a prime of good reduction,
a set $X_{\text{known}} \subseteq X(\Q)$ and a set $X_{\text{extra}} \subseteq X(\Q_p)$ known
to some finite $p$-adic precision, distinct from any of the $X_{\text{known}}$ to that precision.
The goal of the Mordell--Weil sieve, which we describe in this section,
is to describe extra conditions that the points of $X(\Q)$ satisfy that the points in
$X_{\text{extra}}$ do not.
See also \cite{Siksek:mw, BruinStoll,  BBM17}.

We will show that any rational point must be sufficiently close $p$-adically to an element of $X_{\text{known}}$.
To do this, one proves that for each $x \in X(\Q)$, there is some $y\in X_{\text{known}}$
such that $[x-y]\in J(\Q)$ is $p$-adically close to the identity in $J(\Q)$.  We can get a handle on being $p$-adically close to $0 \in J(\Q)$ using the $p$-adic
filtration of $J(\Q_p)$ by 
\[
J_i = \left\{x \in J(\Q_p) \colon x \equiv 0 \pmod{p^i} \right\}\text.
\]
The important property of this filtration that we will make use of is that
\[
J_0 / J_1 \simeq J(\F_p),\; J_{i}/J_{i+1} \simeq \F_p^{\dim J}\,,
\]
so that $p$-adically close rational points must have difference in the Jacobian divisible by a large power of $p$.
Then for any $D\in J(\Q)$ we have $\#J(\F_p)\cdot p^i \cdot D \in J_{i+1}$.

The Mordell--Weil sieve locates small cosets within $J(\Q)$ (that is,
cosets of large index), that contain the image of $X(\Q)$ under the Abel-Jacobi map
$i_b\colon X \to J$ sending $x$ to $[x-b]$.
The sieve plays off local information at a finite set of primes $v$ against the global
Mordell--Weil group structure to find restrictions on $i_b(X(\Q))$.  First we fix a prime $v$ of good reduction and consider the following commutative diagram:
\begin{equation}\label{eqn:diagram-mw-sieve}
\begin{tikzpicture}[->,>=stealth',baseline=(current  bounding  box.center)]
 \node[] (X) {$X(\Q)$};
 \node[right of=X, node distance=3.2cm]  (J) {$J(\Q)$};
  \node[below of=X, node distance=1.5cm]  (Xv) {$X(\F_v)$};
  \node[right of=Xv,node distance=3.2cm] (Jv) {$J(\F_v)$};

 \path (X)  edge node[left]{$\operatorname{red}_{X,v}$} (Xv);
 \path (J)  edge node[right]{$\operatorname{red}_{J,v}$} (Jv);
 \path (X)  edge node[above]{$i_b$} (J);
  \path (Xv)  edge node[below]{$i_{b,v}$} (Jv);
\end{tikzpicture}
\end{equation}

The commutativity of the diagram implies
that the image of $X(\Q)$ along $\operatorname{red}_{J,v}\circ i_b$ is contained in the
image of $i_{b,v}$.
The advantage of this observation is that the bottom row of the diagram deals with finite objects and information about these may be computed effectively.
In particular we can find $\operatorname{im} i_{b,v}$ given equations for $X$. In our
setting of a hyperelliptic curve, algorithms for this go back to~\cite{Cantor87}, and in
general one can make use of work of
Khuri-Makdisi~\cite{KM07}.
Pulling the computed image $\operatorname{im} i_{b,v}$
back to $J(\Q)$ gives a union of cosets for the kernel of
$\operatorname{red}_{J,v}$ that contains the image of $X(\Q)$. We will want to pick $v$ so
that the kernel of this map provides non-trivial information about cosets of the target
subgroup, which means that the index of the kernel is divisible by $p$.
The Mordell--Weil sieve diagram can be amended by using several primes $v$ of good
reduction or working with residue classes of $J(\Q)$; it is also possible to make use of
primes of bad reduction and to go deeper into the filtration $(J_i)_i$. 

For simplicity, we suppose that $r=g$ and we fix a basis $D_1,\ldots, D_g$ of
$J(\Q)/J(\Q)_{\text{torsion}}$.  If $x \in X(\Q_p)$ were to be rational, and we expressed
\[ i_b(x) = \sum_{j = 1}^g m_j D_j,\,m_j\in \Z\,, \]
then we would have, via the linearity of the Coleman integral of regular $1$-forms on the Jacobian,
\begin{equation}\label{eqn:m_j}
 \int_b^x \omega_i = \sum_{j=1}^g m_j \int_0^{D_j} \omega_ i,\,\text{ for each }i\in
\{1,\ldots,g\}\, 
\end{equation}
where we identify $\omega_i$ with the holomorphic differential it induces on $J$ via $\iota_b$.
This can be used to
determine the $m_j$ for given $x\in X(\Q_p)$ modulo $p^n$ for any $n$.
We are done if we can show for every $x\in X_{\text{extra}}$ that the resulting coset
of $J(\Q)/p^nJ(\Q)$ does not meet
the pullback of $i_{b,v}$ under $\operatorname{red}_{J,v}$ for some $v$.

\subsubsection{$X_0(67)^+$}
We now give some details of this computation for $X_0(67)^+$, 
using the model
$$  y^2 = x^6 + 2x^5 + x^4 - 2x^3 + 2x^2 - 4x + 1\,;$$
we have for $X_{\text{known}}$ the $10$ points found in~\eqref{eqn:knownpts67}.
The quadratic Chabauty computation described above also results in a set
$X_{\text{extra}}$ of $11$-adic points of cardinality~14, known to finite precision, whose elements are roots of the
function $Q$ in~\eqref{eq:final_qc_fn}, but which do not appear to be rational. See Table \ref{table:Qp-pointsX067+} for their $x$-coordinates. 

As above, we take $b=(1,1)$; with this choice $D_1 =  i_b(\infty^-)$ and $D_2= i_b(\infty^+)$ are generators
for $J(\Q)$. 
In terms of this basis, we find that $i_b(X_{\text{known}})$ is given by pairs 
\[ (m_1, m_2) \,\in\, \{ ( 1, 0 ), ( 0, 1 ), ( -6, 4 ), ( 7, -3 ), ( 3, -1 ), ( -2, 2 ), (
1, 1 ), ( 0, 0 ), ( 8, -5 ), ( -7, 6 ) \}\,.\]
Since we are working with $p=11$, 
we look for primes $v$ such that $\operatorname{ord}_{11}( \#J(\F_v))$ is large.

We find that
\[J(\F_{31}) \simeq (\Z/(3\cdot 11))^2 \text{ and } J(\F_{137}) \simeq \Z/3 \oplus \Z/(3\cdot 11^2\cdot 19)\]

and  the image of $J(\Q)/11^2J(\Q)$ inside these groups surjects onto the $11$-parts. 
We pull back the images of $i_{b,31}$ and $i_{b,137}$ to 
cosets for $J(\Q)/11^2J(\Q)$. 
Using~\eqref{eqn:m_j} we compute $i_b(x)$ modulo $11^2$ for all $x\in X_{\text{extra}}$,
assuming $x$ is rational, and we find that this does not meet our 
cosets for 31 or 137. 

\subsubsection{Further examples}

In the case of $N=73$ we run computations analogous to the ones described above, using the prime
$p =37$ for  the quadratic Chabauty procedure, and 
applying the Mordell--Weil sieve with the prime $9511$ to rule out the extra $37$-adic points.
Likewise for $N = 103$ we can perform
quadratic Chabauty at $p =3$.
This gives  6 ``extra'' $3$-adic points, which can be shown
to be non-rational by applying the Mordell--Weil sieve using the prime $397$.

\subsubsection{Conclusion}
In summary, we have shown:
\begin{theorem}
    The number of rational points on the Atkin--Lehner quotient modular curves $X_0(N)^+$
    for $N\in \{67, 73, 103\}$ are as follows:
    \[
        \#X_0(67)^+(\Q) = 10\,,\quad 
        \#X_0(73)^+(\Q) = 10\,,\quad
        \#X_0(103)^+(\Q) = 8 \,.
    \]
\end{theorem}

According to~\cite{Gal96}, this shows that $X_0(67)^+(\Q)$ contains no exceptional points and that
$X_0(73)^+(\Q)$ and $X_0(103)^+(\Q)$  contain precisely one exceptional point each, up to the hyperelliptic
involution. Here an exceptional point is a rational point that is neither a cusp nor a
CM point.

Furthermore, we may conclude that the table in \cite[\S~4.6]{Box19} contains all quadratic
points on $X_0(67)$ and the table in \cite[\S~4.7]{Box19} contains all quadratic
points on $X_0(73)$, complementing \cite[Theorem~1.1]{Box19}.

Finally, our theorem implies that the list of $j$-invariants of $\Q$-curves attached to non-cuspidal rational
points on $X_0(N)^+$ given in~\cite[\S4.1]{BGX19} is complete for $N\in \{67,73,103\}$.

\subsection*{Acknowledgements}
We would like to thank Barry Mazur, whose kind enthusiasm and encouragement were a great
motivation for the writing of this paper. We thank Levent Alpoge, Netan Dogra, Kiran Kedlaya, Minhyong
Kim, Bjorn Poonen, Michael Stoll, Akshay Venkatesh, and John Voight for many useful conversations and suggestions relating to the material in this article, and Bas Edixhoven and Guido Lido for sharing an early version of their preprint~\cite{EL19} with us and for providing several useful suggestions and corrections on an earlier version of this article. Most of the computations described in Section~\ref{sec:Examples} were done at the Workshop ``Arithmetic Statistics and Diophantine Stability'' at the Fondation des Treilles in July 2018. We are grateful to Barry Mazur and Karl Rubin for organizing this event.

JB was supported in part by NSF grant DMS-1702196, the Clare Boothe Luce Professorship (Henry Luce Foundation), and
Simons Foundation grant \#550023. 
AB was supported by Simons Foundation grant \#550023.
FB was supported by EPSRC and by Balliol College through a Balliol Dervorguilla scholarship.
BL was supported by an NSF postdoctoral fellowship. 
SM was partially supported by DFG grant MU 4110/1-1 and an NWO Vidi grant.
NT was supported by an NSF Graduate Research Fellowship under grant \#1122374, by
Simons Foundation grant \#550033, and by the  Research and Training Group in Algebra, Algebraic Geometry, and Number Theory at the University of Georgia during various stages of this project.
JV was supported by Francis Brown and ERC-COG 724638 `GALOP', the
Carolyn and Franco Gianturco Fellowship at Linacre College (Oxford),
and NSF Grant No. DMS-1638352, during various stages of this project.
\bibliographystyle{alpha}

\bibliography{References}

\newcommand{\etalchar}[1]{$^{#1}$}
\begin{thebibliography}{BDCKW18}

\bibitem[Bal15]{B15}
J.~Balakrishnan.
\newblock {C}oleman integration for even-degree models of hyperelliptic curves.
\newblock {\em LMS J. Comput. Math.}, 18(1):258--265, 2015.

\bibitem[BBBM19]{BBBM}
J.S. Balakrishnan, A.~Besser, F.~Bianchi, and J.S. M\"uller.
\newblock Explicit quadratic {C}habauty over number fields.
\newblock {\em ArXiv preprint}, arXiv:1910.04653, 2019.

\bibitem[BBK10]{BBK10}
J.~Balakrishnan, R.~Bradshaw, and K.~Kedlaya.
\newblock Explicit {C}oleman integration for hyperelliptic curves.
\newblock In {\em Algorithmic number theory (ANTS-IX)}, volume 6197 of {\em
  Lecture Notes in Comput. Sci.}, pages 16--31. Springer, Berlin, 2010.

\bibitem[BBM17]{BBM17}
Jennifer~S. Balakrishnan, Amnon Besser, and J.~Steffen M\"{u}ller.
\newblock Computing integral points on hyperelliptic curves using quadratic
  {C}habauty.
\newblock {\em Math. Comp.}, 86(305):1403--1434, 2017.

\bibitem[BCP97]{BCP97}
W.~Bosma, J.~Cannon, and C.~Playoust.
\newblock {The Magma algebra system I: The user language}.
\newblock {\em J. Symb. Comp}, 24(3-4):235--265, 1997.

\bibitem[BD18]{BD18}
J.~Balakrishnan and N.~Dogra.
\newblock Quadratic {C}habauty and rational points {I}: $p$-adic heights.
\newblock {\em Duke Math. J.}, 167(11):1981--2038, 2018.

\bibitem[BD19a]{BD19}
J.~Balakrishnan and N.~Dogra.
\newblock Quadratic {C}habauty and rational points {II}: Generalised height
  functions on {S}elmer varieties.
\newblock {\em Arxiv preprint}, arXiv:1705.00401v2, 2019.

\bibitem[BD19b]{BeD19}
A.~Betts and N.~Dogra.
\newblock Ramification of {\'e}tale path torsors and harmonic analysis on
  graphs.
\newblock {\em ArXiv preprint}, arXiv:1909.05734, 2019.

\bibitem[BDCKW18]{BDCKW18}
Jennifer~S. Balakrishnan, Ishai Dan-Cohen, Minhyong Kim, and Stefan Wewers.
\newblock A non-abelian conjecture of {T}ate-{S}hafarevich type for hyperbolic
  curves.
\newblock {\em Math. Ann.}, 372(1-2):369--428, 2018.

\bibitem[BDM{\etalchar{+}}19]{BDMTV19}
J.~Balakrishnan, N.~Dogra, S.~M\"uller, J.~Tuitman, and J.~Vonk.
\newblock Explicit {C}habauty--{K}im for the split {C}artan modular curve of
  level $13$.
\newblock {\em Annals of Math.}, 189(3), 2019.

\bibitem[BE19]{BE19}
D.~Bertrand and B.~Edixhoven.
\newblock Pink's conjecture on unlikely intersections and families of
  semi-abelian varieties.
\newblock {\em ArXiv Preprint}, arXiv:1904.01788, 2019.

\bibitem[BGX19]{BGX19}
Francesc Bars, Josep González, and Xavier Xarles.
\newblock Hyperelliptic parametrizations of $\mathbb{Q}$-curves.
\newblock {\em ArXiv preprint}, arXiv:1903.05032v2, 2019.

\bibitem[Box19]{Box19}
J.~Box.
\newblock Quadratic points on modular curves with infinite {M}ordell--{W}eil
  group.
\newblock {\em ArXiv Preprint}, arXiv:1906.05206, 2019.

\bibitem[BPSS]{g3hyp}
A.~Booker, D.~Platt, J.~Sijsling, and A.~Sutherland.
\newblock Genus 3 hyperelliptic curves.
\newblock \url{http://math.mit.edu/~drew/lmfdb\_genus3\_hyperelliptic.txt}.

\bibitem[BS10]{BruinStoll}
Nils Bruin and Michael Stoll.
\newblock The {M}ordell-{W}eil sieve: proving non-existence of rational points
  on curves.
\newblock {\em LMS J. Comput. Math.}, 13:272--306, 2010.

\bibitem[BT17]{BT18}
Benjamin Bakker and Jacob Tsimerman.
\newblock The ax--schanuel conjecture for variations of hodge structures.
\newblock {\em ArXiv preprint}, arXiv:1712.05088, 2017.

\bibitem[BT19]{BT19}
J.~Balakrishnan and J.~Tuitman.
\newblock Explicit {C}oleman integration for curves.
\newblock {\em Preprint}, 2019.

\bibitem[Can87]{Cantor87}
David~G. Cantor.
\newblock Computing in the {J}acobian of a hyperelliptic curve.
\newblock {\em Math. Comp.}, 48(177):95--101, 1987.

\bibitem[CG89]{CG89}
Robert~F. Coleman and Benedict~H. Gross.
\newblock {$p$}-adic heights on curves.
\newblock In {\em Algebraic number theory}, volume~17 of {\em Adv. Stud. Pure
  Math.}, pages 73--81. Academic Press, Boston, MA, 1989.

\bibitem[Cha41]{Cha41}
C.~Chabauty.
\newblock Sur les points rationels des courbes alg\'ebriques de genre
  sup\'erieur \`a l'unit\'e.
\newblock {\em C.R. Acad. Sci.}, 212:882--884, 1941.

\bibitem[CK10]{CK10}
J.~Coates and M.~Kim.
\newblock Selmer varieties for curves with {CM} jacobians.
\newblock {\em Kyoto J. Math.}, 50(4):827--852, 2010.

\bibitem[CMSV19]{CMSV}
Edgar Costa, Nicolas Mascot, Jeroen Sijsling, and John Voight.
\newblock Rigorous computation of the endomorphism ring of a {J}acobian.
\newblock {\em Math. Comp.}, 88(317):1303--1339, 2019.

\bibitem[Coh00]{Coh00}
H.~Cohen.
\newblock {\em Advanced Topics in Computational Number Theory}, volume 193 of
  {\em Graduate Texts in Mathematics}.
\newblock Springer, 2000.

\bibitem[Col85]{Col85}
R.~Coleman.
\newblock Torsion points on curves and $p$-adic abelian integrals.
\newblock {\em Annals of Math.}, 121:111--168, 1985.

\bibitem[Cre97]{Cre97}
J.~Cremona.
\newblock {\em Algorithms for modular elliptic curves}.
\newblock Cambridge University Press, Cambridge, second edition, 1997.

\bibitem[CS86]{CS86}
G.~Cornell and J.~Silverman.
\newblock {\em Arithmetic Geometry}.
\newblock Springer-Verlag, New York, 1986.

\bibitem[DC17]{DC}
Ishai Dan-Cohen.
\newblock Mixed {T}ate motives and the unit equation {II}.
\newblock {\em Arxiv preprint}, 2017.

\bibitem[DCW15]{DCW15}
Ishai Dan-Cohen and Stefan Wewers.
\newblock Explicit {C}habauty-{K}im theory for the thrice punctured line in
  depth 2.
\newblock {\em Proc. Lond. Math. Soc. (3)}, 110(1):133--171, 2015.

\bibitem[DCW16]{DCW16}
Ishai Dan-Cohen and Stefan Wewers.
\newblock Mixed {T}ate motives and the unit equation.
\newblock {\em Int. Math. Res. Not. IMRN}, (17):5291--5354, 2016.

\bibitem[Del85]{Deligne_ladic}
Pierre Deligne.
\newblock Repr\'esentations $\ell$-adiques.
\newblock In {\em Ast\'erisque 124}. 1985.

\bibitem[Del89]{Del89}
P.~Deligne.
\newblock Le groupe fondamental de la droite projective moins trois points.
\newblock In {\em Galois groups over {$\mathbf{Q}$}}, volume~16 of {\em Math.
  Inst. Res. Inst. Publ.}, pages 79--297. Springer-Verlag, 1989.

\bibitem[Dog19]{D19}
N.~Dogra.
\newblock Unlikely intersections and the {C}habauty--{K}im method over number
  fields.
\newblock {\em ArXiv preprint}, arXiv:1903.05032v2, 2019.

\bibitem[Dok04]{Dok04}
Tim Dokchitser.
\newblock Computing special values of motivic {$L$}-functions.
\newblock {\em Experiment. Math.}, 13(2):137--149, 2004.

\bibitem[DR73]{DR73}
P.~Deligne and M.~Rapoport.
\newblock Les sch\'emas de modules de courbes elliptiques.
\newblock In W.~Kuyk, editor, {\em Modular forms in one variable II}, volume
  349 of {\em LNM}, pages 143--316. Springer-Verlag, 1973.

\bibitem[EH17]{EH}
J.~Ellenberg and D.~Hast.
\newblock Rational points on solvable curves over {$\mathbb{Q}$} via
  non-abelian {C}habauty.
\newblock {\em ArXiv preprint}, arXiv:1706.00525, 2017.

\bibitem[EL19]{EL19}
B.~Edixhoven and G.~Lido.
\newblock Geometric quadratic {C}habauty.
\newblock {\em ArXiv Preprint}, arXiv:1910.10752, 2019.

\bibitem[Fal83]{Fal83}
G.~Faltings.
\newblock Endlichkeitss\"atze f\"ur abelsche {V}ariet\"aten \"uber
  {Z}ahlk\"orpern.
\newblock {\em Invent. Math.}, 73(3):349--366, 1983.

\bibitem[FL82]{FL82}
Jean-Marc Fontaine and Guy Laffaille.
\newblock Construction de repr\'esentations $p$-adiques.
\newblock {\em Annales scientifiques de l'\'E.N.S.}, 14(4):547--608, 1982.

\bibitem[FM12]{FM12}
Benson Farb and Dan Margalit.
\newblock {\em A Primer on Mapping Class Groups}.
\newblock Princeton University Press, 2012.

\bibitem[Gal96]{Gal96}
S.\thinspace{}D. Galbraith.
\newblock Equations for modular curves.
\newblock {\em Oxford DPhil thesis}, 1996.

\bibitem[GLLM15]{GLLM}
Fritz Grunewald, Michael Larsen, Alexander Lubotzky, and Justin Malestein.
\newblock Arithmetic quotients of the mapping class group.
\newblock {\em Geometric and Functional Analysis}, 25:1493--1542, 2015.

\bibitem[Gro97]{Gro97}
Alexander Grothendieck.
\newblock Brief an {G.} {F}altings.
\newblock In {\em Geometric Galois actions, 1}, volume 242 of {\em London Math.
  Soc. Lecture Note Ser.}, pages 49--58. Cambridge University Press, 1997.

\bibitem[GZ86]{GZ86}
B.~Gross and D.~Zagier.
\newblock Heegner points and derivatives of {L}-series.
\newblock {\em Invent. Math.}, 84(2):225--320, 1986.

\bibitem[Had11]{Had11}
M.~Hadian.
\newblock Motivic fundamental groups and integral points.
\newblock {\em Duke Math. J.}, 160:503--565, 2011.

\bibitem[Har12]{harrison}
M.\thinspace{}C. Harrison.
\newblock An extension of {K}edlaya{'}s algorithm for hyperelliptic curves.
\newblock {\em J. Symb. Comp.}, 47(1):89 -- 101, 2012.

\bibitem[Ked01]{Ked01}
Kiran~S. Kedlaya.
\newblock Counting points on hyperelliptic curves using {M}onsky-{W}ashnitzer
  cohomology.
\newblock {\em J. Ramanujan Math. Soc.}, 16(4):323--338, 2001.

\bibitem[Ked07]{Ked07}
K.~Kedlaya.
\newblock $p$-{Adic} cohomology: from theory to practice.
\newblock {\em Arizona Winter School Notes}, 2007.

\bibitem[Ked09]{Ked05}
Kiran~S. Kedlaya.
\newblock {$p$}-adic cohomology.
\newblock In {\em Algebraic geometry---{S}eattle 2005. {P}art 2}, volume~80 of
  {\em Proc. Sympos. Pure Math.}, pages 667--684. Amer. Math. Soc., Providence,
  RI, 2009.

\bibitem[Kim05]{Kim05}
M.~Kim.
\newblock The motivic fundamental group of $\mathbf{P}^1 \backslash \{
  0,1,\infty \}$ and the theorem of {S}iegel.
\newblock {\em Invent. Math.}, 161:629--656, 2005.

\bibitem[Kim09]{Kim09}
M.~Kim.
\newblock The unipotent {A}lbanese map and {S}elmer varieties for curves.
\newblock {\em Publ. RIMS}, 45:89--133, 2009.

\bibitem[Kim10]{Kim10b}
M.~Kim.
\newblock Massey products for elliptic curves of rank $1$.
\newblock {\em J. Amer. Math. Soc.}, 23(3):725--747, 2010.

\bibitem[KL89]{KL89}
V.~A. Kolyvagin and D.~Yu. Logach{\"e}v.
\newblock Finiteness of the {S}hafarevich-{T}ate group and the group of
  rational points for some modular abelian varieties.
\newblock {\em Algebra i Analiz}, 1(5):171--196, 1989.

\bibitem[KM04]{KM03}
Kamal Khuri-Makdisi.
\newblock Linear algebra algorithms for divisors on an algebraic curve.
\newblock {\em Math.\ Comp.}, 73:333--357, 2004.

\bibitem[KM07]{KM07}
Kamal Khuri-Makdisi.
\newblock Asymptotically fast group operations on {J}acobians of general
  curves.
\newblock {\em Math. Comp.}, 76(260):2213--2239, 2007.

\bibitem[Lau04]{Lau04}
Alan G.~B. Lauder.
\newblock Deformation theory and the computation of zeta functions.
\newblock {\em Proc. London Math. Soc. (3)}, 88(3):565--602, 2004.

\bibitem[Lau06]{Lau06}
Alan G.~B. Lauder.
\newblock A recursive method for computing zeta functions of varieties.
\newblock {\em LMS J. Comput. Math.}, 9:222--269, 2006.

\bibitem[Loo97]{Looijenga}
Eduard Looijenga.
\newblock Prym representations of mapping class groups.
\newblock {\em Geom. Dedicata}, 64(1):69--83, 1997.

\bibitem[LV18]{LV18}
B.~Lawrence and A.~Venkatesh.
\newblock Diophantine problems and $p$-adic period mappings.
\newblock {\em ArXiv preprint}, arXiv:1807.02721v1, 2018.

\bibitem[Mas19]{Mascot}
Nicolas Mascot.
\newblock Hensel-lifting torsion points on {J}acobians and {G}alois
  representations.
\newblock {\em ArXiv preprint}, arXiv:1808.03939, 2019.

\bibitem[Mor22]{Mor22}
L.~J. Mordell.
\newblock On the rational solutions of the indeterminate equations of the third
  and fourth degrees.
\newblock {\em Proc. Cambridge Phil. Soc.}, 21:179--192, 1922.

\bibitem[MP12]{McP12}
W.~McCallum and B.~Poonen.
\newblock The method of {C}habauty and {C}oleman.
\newblock In {\em Explicit methods in number theory}, volume~36 of {\em Panor.
  Synth\`eses}, pages 99--117. Soc. Math. France, Paris, 2012.

\bibitem[Mum83]{Tata1}
David Mumford.
\newblock {\em Tata lectures on theta}.
\newblock Birkh\"auser, 1983.

\bibitem[Mur92]{Murabayashi}
Naoki Murabayashi.
\newblock On normal forms of modular curves of genus {$2$}.
\newblock {\em Osaka J. Math.}, 29(2):405--418, 1992.

\bibitem[Nek93]{Nek93}
J.~Nekovar.
\newblock On $p$-adic height pairings.
\newblock In {\em S\'eminaire de {T}h\'eorie des {N}ombres, {P}aris 1990-1991},
  pages 127--202. Birkh\"auser, 1993.

\bibitem[Ogg75]{Ogg75}
A.~P. Ogg.
\newblock Automorphismes de courbes modulaires.
\newblock In {\em S\'{e}minaire {D}elange-{P}isot{P}oitou (16e ann\'{e}e:
  1974/75), {T}h\'{e}orie des nombres, {F}asc. 1, {E}xp. {N}o. 7}, page~8.
  1975.

\bibitem[Sik13]{siksek:nf}
Samir Siksek.
\newblock Explicit {C}habauty over number fields.
\newblock {\em Algebra Number Theory}, 7(4):765--793, 2013.

\bibitem[Sik15]{Siksek:mw}
Samir Siksek.
\newblock Chabauty and the {M}ordell-{W}eil sieve.
\newblock In {\em Advances on superelliptic curves and their applications},
  volume~41 of {\em NATO Sci. Peace Secur. Ser. D Inf. Commun. Secur.}, pages
  194--224. IOS, Amsterdam, 2015.

\bibitem[ST18]{ST}
Nick Salter and Bena Tshishiku.
\newblock Arithmeticity of the monodromy of some {K}odaira fibrations.
\newblock {\em ArXiv Preprint}, arXiv:1805.06789, 2018.

\bibitem[Ste00]{Ste00}
W.~Stein.
\newblock {\em Explicit approaches to modular abelian varieties}.
\newblock ProQuest LLC, Ann Arbor, MI, 2000.
\newblock Thesis (Ph.D.)--University of California, Berkeley.

\bibitem[SV14]{SV14}
J.~Sijsling and J.~Voight.
\newblock On computing {B}ely\u{\i} maps.
\newblock {\em Publications math\'ematiques de Besan\c{c}on}, (1):73--131,
  2014.

\bibitem[Tui16]{Tui16}
Jan Tuitman.
\newblock Counting points on curves using a map to {$\mathbf{P}^1$}.
\newblock {\em Math. Comp.}, 85(298):961--981, 2016.

\bibitem[Tui17]{Tui17}
J.~Tuitman.
\newblock Counting points on curves using a map to {$\bold{P}^1$}, {II}.
\newblock {\em Finite Fields Appl.}, 45:301--322, 2017.

\end{thebibliography}


\end{document}